\documentclass[twocolumn]{autart}    

\usepackage{graphicx}          

\usepackage{amsmath}
\usepackage{mathrsfs}
\usepackage{amssymb}
\usepackage{booktabs}
\usepackage{tabularx}
\usepackage{algorithm} 
\usepackage{algorithmic}
\makeatletter
\@addtoreset{thm}{section}

\makeatother

\begin{document}

\begin{frontmatter}

\title{Sequential feedback optimization with application to wind farm control} 

\thanks[footnoteinfo]{This work is supported by the Dutch Organization for Scientific Research (NWO) under the project "Online Optimization for Offshore Wind Farms" (OCENW.M.21.190). Corresponding author S.~Grammatico.}

\author[Delft]{Shijie Huang}\ead{S.Huang-5@tudelft.nl},    
\author[Delft]{Sergio Grammatico}\ead{S.Grammatico@tudelft.nl}             

\address[Delft]{Delft Center for Systems and Control, TU Delft}  


\begin{abstract}                          
This paper develops a sequential-linearization feedback optimization framework for driving nonlinear dynamical systems to an optimal steady state. A fundamental challenge in feedback optimization is the requirement of accurate first-order information of the steady-state input-output mapping, which is computationally prohibitive for high-dimensional nonlinear systems and often leads to poor performance when approximated around a fixed operating point. To address this limitation, we propose a sequential  algorithm that adaptively updates the linearization point during optimization, maintaining local accuracy throughout the trajectory. We prove convergence to a neighborhood of the optimal steady state with explicit error bounds. To reduce the computational burden of repeated linearization operations, we further develop a multi-timescale variant where linearization updates occur at a slower timescale than optimization iterations, achieving significant computational savings while preserving convergence guarantees. The effectiveness of the proposed framework is demonstrated via numerical simulations of a realistic wind farm control problem. The results validate both the theoretical convergence predictions and the expected computational advantages of our multi-timescale formulation.
\end{abstract}

\end{frontmatter}

\section{Introduction}
Operating complex dynamical systems at their optimal steady state is a fundamental challenge across various engineering applications, from chemical process control to power grid management and renewable energy systems \cite{francois2013measurement}, \cite{krishnamoorthy2022real}, \cite{kebir2019real}. The goal is to find control inputs that drive the system to an equilibrium point optimizing performance metrics while satisfying operational constraints. This steady-state optimization problem becomes particularly challenging for nonlinear, high-dimensional systems subject to unknown disturbances, where traditional feedforward optimization approaches require accurate knowledge of the steady-state input-output mapping and disturbances --- information that is often computationally intractable or practically unavailable in real applications \cite{he2023model}, \cite{molzahn2017survey}.

These challenges are particularly evident in wind farm control. As wind power becomes increasingly critical in global sustainable energy systems, modern wind farms must coordinate dozens to even hundreds of turbines to maximize steady collective power output while managing complex aerodynamic wake interactions \cite{boersma2017tutorial}. Most existing optimization approaches use low-fidelity steady-state flow model such as FLORIS \cite{gebraad2016wind} to enable tractable solutions, but suffer from model-reality mismatches, which might lead to suboptimal performance when tested in high-fidelity models \cite{kheirabadi2019quantitative}. Higher-fidelity flow models like WFSim \cite{boersma2018control}, governed by the Navier--Stokes equations with typical high-dimensional state, can better capture wake dynamics but hinder the direct application of traditional optimization methods. Economic model predictive framework is usually adopted in the literature \cite{vali2019adjoint}, but it faces computational inefficiencies and lack theoretical convergence to the optimal steady state. 

Feedback optimization (FO) offers a promising alternative framework to address this issue \cite{hauswirth2024optimization}. By integrating optimization algorithms directly with real-time measurements, FO methods can drive systems to optimal steady states without requiring perfect model knowledge. The approach replaces model-based steady-state predictions with actual system measurements, enabling optimization algorithms to adapt to real system behavior. This framework has been extensively studied theoretically and demonstrated success in applications such as frequency control \cite{kasis2017stability} and voltage control \cite{dominguez2023online} in power systems, establishing FO as a powerful tool for steady-state optimization of complex systems.

However, a fundamental limitation of existing FO algorithms is their dependence on the Jacobian matrix of the input-output mapping. While the Jacobian matrix is readily available for linear systems, its computation becomes computationally prohibitive for complex nonlinear dynamical systems, particularly those with high-dimensional state spaces. Previous efforts to address this challenge have focused on approximate feedback optimization approaches that employ surrogate Jacobian matrices obtained through linearization around fixed operating points \cite{colombino2019towards}, \cite{ortmann2020experimental}, \cite{colot2024optimal}. While this approach enables a simple and tractable implementation, the quality of the resulting approximate solution is highly dependent on how well the chosen linearization point aligns with the unknown global optimal steady state. In wind farm control applications, where optimal operating conditions depend on complex, spatially-distributed aerodynamic interactions, poor linearization point selection can cause algorithms to converge to suboptimal solutions that significantly underperform the global optimum.
 
 This paper addresses this fundamental limitation by developing a novel sequential feedback optimization framework to eliminate the dependence on a-priori operating point selection, while maintaining theoretical convergence guarantees.

\subsection{Related work}
\textbf{Feedback Optimization.} The key idea of feedback optimization is to implement the traditional optimization algorithms in a feedback loop with the original dynamical system by using real-time measurements of the output to replace the steady output, thus eliminating the requirement on perfect steady-state input-output mapping. \cite{menta2018stability} analyzed the stability of feedback gradient descent flows and apply it to power systems. \cite{lawrence2018optimal} extended these concepts by combining proportional-integral control with feedback gradient methods for linear time-invariant systems. The framework has been generalized to address time-varying optimization problems \cite{bernstein2019online}, \cite{colombino2019online}, where the authors established tracking properties for optimal trajectories under changing conditions. Recent advances have tackled nonlinear dynamical systems \cite{chen2023online}, \cite{hauswirth2020timescale}, \cite{he2024online}, providing stability guarantees under various assumptions on system dynamics. \cite{hauswirth2024optimization} provided a comprehensive survey of these developments.

 A persistent challenge in feedback optimization is the requirement for the Jacobian matrix of the steady-state input-output mapping. \cite{colombino2019towards} addresses this by proposing linearization around fixed operating points and analyzes the robustness of such approximations, with subsequent experimental validation by \cite{ortmann2020experimental}. Furthermore, \cite{picallo2022adaptive} takes a data-driven approach using recursive least-squares to estimate sensitivities online, though requires learning the sensitivity from noisy measurements over multiple iterations and careful tuning of stochastic noise parameters.
 
 \textbf{Wind Farm Control.} The effectiveness of wind farm control strategies fundamentally depends on the model fidelities. Low-fidelity steady-state wake models, including the Jensen model \cite{jensen1983note}, Park model \cite{katic1987simple}, and the more recent FLORIS framework \cite{gebraad2016wind}, enable computationally efficient optimization. These models have been leveraged for various optimization algorithms: \cite{park2015cooperative} developed sequential convex programming methods, while \cite{annoni2019efficient} proposed distributed optimization using proximal primal-dual algorithms. However, high-fidelity simulation and wind tunnel testing have shown that these simplified models can lead to suboptimal performance due to unmodeled dynamics \cite{annoni2016analysis}, \cite{campagnolo2016wind}.
 
Medium-fidelity dynamic models bridge the gap between computational efficiency and physical accuracy. The Wind Farm Simulator (WFSim) \cite{boersma2018control}, based on two-dimensional Navier--Stokes equations, captures wake dynamics while remaining tractable for control design. \cite{vali2019adjoint} demonstrates adjoint-based model predictive control using WFSim. High-fidelity models based on large-eddy simulations provide the most accurate representations but remain limited to offline validation \cite{bay2019unlocking}, \cite{andersen2020global} and parameter identification \cite{gebraad2016wind} due to computational constraints. Beyond the model-based methods mentioned above, model-free approaches such as extremum seeking control have also been applied to wind farm \cite{johnson2012assessment}, \cite{ciri2017model}. For a more comprehensive review for various wind farm control strategies, we refer the interested reader to \cite{andersson2021wind}. 

\subsection{Contributions}
The main contributions of this paper are threefold:

\textbf{Sequential Adaptive Linearization and Convergence Analysis:} We propose a novel sequential feedback optimization algorithm that adaptively updates the linearization point at each iteration and establish rigorous theoretical convergence guarantees. Specifically, we prove that despite using time-varying linearizations, the state-control sequence converges to a neighborhood of the optimal steady state. We derive explicit bounds on the steady-state approximation error as a function of system parameters.

\textbf{Multi-Timescale Computational Framework:} To address the computational burden in high-dimensional systems, we develop a multi-timescale variant that performs linearization updates at a slower timescale, significantly reducing computational overhead while maintaining convergence properties with theoretical analysis of the trade-off between efficiency and accuracy.

\textbf{Numerical Application to Wind Farm Control:} We demonstrate the practical effectiveness of our approach through extensive simulations using a realistic offshore wind farm configuration, showing substantial improvements in steady-state power production compared to the greedy control strategy while validating the theoretical predictions regarding convergence behavior and computational trade-offs.

This paper is structured as follows. Section \ref{section:problem} introduces the problem setup. Section \ref{section:algorithm} presents algorithm design and main theoretical results, while Section \ref{section:multi-scale} further provides the multi-timescale variant. Section \ref{section:application} formulates
the wind farm power maximization problem and gives both medium-fidelity and high-fidelity
numerical simulations. Finally, Section \ref{section:conclusion} concludes the paper.

\textbf{Notations.} We denote $\mathbb{R}^n$ as the $n$-dimensional real Euclidean space. For a column vector $x\in\mathbb{R}^n$ (matrix $A\in \mathbb{R}^{m\times n}$, $x^{\top}$ ($A^{\top}$) denotes its transpose. Denote $\text{proj}_{\mathcal{U}}$ as the Euclidean projection operator on a set $\mathcal{U}$, i.e., $\text{proj}_{\mathcal{U}}(v):= \arg\min_{u\in\mathcal{U}}\{\|u-v\|\}$. For a multi-variable function $f(x)$, denote $\nabla f(x)$ as the gradient.

\section{Problem formulation}\label{section:problem}
We consider the problem to find an efficient steady state of a plant with nonlinear discrete-time dynamics:
\begin{equation}\label{plant_model}
\left\{\begin{array}{rcl}
x^{+} &=& f(x,u) + w_1\\
y\ \ \  &=& g(x,u) + w_2
\end{array}\right.
\end{equation}
where $x\in\mathbb{R}^n$ represents the system state, $u\in\mathcal{U}\subseteq\mathbb{R}^p$ represents the control variable, $y\in\mathbb{R}^m$ represents the measured output, and $w = \mathrm{col}\{w_1,w_2\}$ represents unknown disturbances. We assume that $\mathcal{U}$ is a convex and compact set representing physical constraints on the control inputs, $f: \mathbb{R}^n\times \mathbb{R}^p\to \mathbb{R}^n$ is a continuously differentiable vector field and the matrix $I - \nabla_x f(x,u)$ is invertible for all $x$ and $u$. The implicit function theorem \cite{khalil2002nonlinear} ensures the existence of a steady-state mapping $\phi: \mathcal{U}\to \mathbb{R}^n$ such that
\[\phi(u,w_1) = f(\phi(u,w_1),u) + w_1,\ \text{for all }\ u\in\mathcal{U}.\]
Furthermore, the steady-state input-output mapping is given by
\begin{equation}\label{input_output}
y_{\text{ss}} = g(\phi(u,w_1),u) + w_2 =: h(u,w).
\end{equation}
This framework is particularly relevant for wind farm control applications under constant inflow conditions, where $x(k)$ represents the spatially discretized wind velocity field, $u(k)$ denotes control inputs including axial induction factors and yaw angles, and $y(k)$ captures performance outputs such as power generation and structural loads \cite{vali2019adjoint}, \cite{boersma2017tutorial}. The nonlinear dynamics $f$ encode the wake propagation physics, typically derived from the Navier--Stokes equations \cite{boersma2018control}.

Formally, we make the following technical assumption on the system in \eqref{plant_model}.
\begin{assum}\label{ass:system}
The system dynamics in \eqref{plant_model} and steady-state mapping in \eqref{input_output} satisfy:
\begin{enumerate}
    \item[(i)] The mapping $x \mapsto f(x,u)$ is uniformly contractive with factor $\rho_f \in (0,1)$ for all $u \in \mathbb{R}^p$, i.e., 
    \[\sup_{(x,u) \in \mathbb{R}^n \times \mathcal{U}}\|\nabla_x f(x,u)\| \leq \rho_f < 1;\]
    
    \item[(ii)] The partial gradients $\nabla_x f(x,u)$ and $\nabla_u f(x,u)$ are uniformly bounded (i.e., $\|\nabla_u f(x,u)\|\le G_u^f$) and Lipschitz continuous with constants $L_{f,x}$ and $L_{f,u}$ respectively, i.e., for all $x_1,x_2 \in\mathbb{R}^n$ and $u_1,u_2 \in\mathcal{U}$:
    \begin{align*}
    \|\nabla_x f(x_1,u_1) - \nabla_x f(x_2,u_2)\|&\leq L_{f,x} (\|x_1 - x_2\|\\
    &\quad + \|u_1 - u_2\|), \\
    \|\nabla_u f(x_1,u_1) - \nabla_u f(x_2,u_2)\|&\leq L_{f,u} (\|x_1 - x_2\|\\
    &\quad + \|u_1 - u_2\|).
    \end{align*}
    
    \item[(iii)] The steady-state input-output mapping $h: \mathcal{U} \to \mathbb{R}^m$ in \eqref{input_output} is $L_h$-Lipschitz continuous:
    \begin{align*}
    \|h(u_1) - h(u_2)\| &\leq L_h \|u_1 - u_2\|, \  \text{for all}\ u_1,u_2 \in \mathcal{U}.
    \end{align*}
\end{enumerate}
\end{assum}
Given a control input $u$, Assumption \ref{ass:system}(i) ensures that the system state converges to the corresponding steady state at an exponential rate, which is essential for the convergence analysis of feedback optimization algorithms \cite[Assumption 2.1]{hauswirth2020timescale}, \cite[Assumption 1]{he2023model}. Conditions (ii)-(iii) provide the regularity needed for the steady-state input-output mapping and its Jacobian, which are fundamental for gradient-based optimization updates \cite[Assumption 1]{chen2023online}.

In this paper, the goal is to optimize over the steady state to ensure that the system settles at a desirable operating point. For instance, in wind farm control under constant inflow conditions, we seek to maximize the steady-state power production or minimize fatigue loads while respecting operational constraints. Specifically, we consider a constrained optimization problem over the steady state $(\bar{u},\bar{y})$ of system \eqref{plant_model}:
\begin{equation}\label{eq:optimization}
\begin{cases}
\min\limits_{\bar{u},\bar{y}} J(\bar{u},\bar{y})\\
\text{s.t.}\ \ \bar{y} = h(\bar{u},w)\\
\quad\ \ \ \bar{u}\in\mathcal{U}.
\end{cases}
\end{equation}
An optimal solution of this optimization problem $(\bar{u}^{\ast},\bar{y}^{\ast})$ is referred to as an optimal steady state-input of system \eqref{plant_model}.

Let us impose the following technical assumption on the cost function.
\begin{assum}\label{ass:cost}
The cost function $J: \mathcal{U} \times \mathbb{R}^m \rightarrow \mathbb{R}$ in \eqref{eq:optimization} satisfies the following conditions:
\begin{enumerate}
    \item[(i)] The partial gradients are uniformly bounded and Lipschitz continuous in $(u,y)$: There exist constants $G_u^J, G_y^J, L_{J,u}, L_{J,y} > 0$ such that  
    \begin{align*}
    &\|\nabla_u J(u_1,y_1)\|\leq G_u^J, \|\nabla_y J(u_1,y_1)\|\leq G_y^J,\\
    &\|\nabla_u J(u_1,y_1) - \nabla_u J(u_2,y_2)\| \leq L_{J,u} (\|u_1 - u_2\|\\
    &\qquad\qquad\qquad\qquad\qquad\qquad\qquad\qquad + \|y_1 - y_2\|),\\
    &\|\nabla_y J(u_1,y_1) - \nabla_y J(u_2,y_2)\| \leq L_{J,y} (\|u_1 - u_2\|\\
    &\qquad\qquad\qquad\qquad\qquad\qquad\qquad\qquad + \|y_1 - y_2\|),
    \end{align*}
    for all $u_1,u_2 \in \mathcal{U}$ and $y_1,y_2 \in \mathbb{R}^m$.
    \item[(ii)] For any fixed $y\in\mathbb{R}^m$, the mapping $u\mapsto\nabla_u J(u,y)$ is $\mu_J$-strongly monotone: 
    \[(\nabla_u J(u_1,y) - \nabla_u J(u_2,y))^{\top}(u_1 - u_2) \geq \mu_J\|u_1 - u_2\|^2,\]
    for all $u_1,u_2 \in \mathcal{U}$ with $\mu_J > 0$.
\end{enumerate}
\end{assum}
Assumption \ref{ass:cost} is commonly employed in the analysis of first-order optimization methods \cite{boyd2004convex} and appears frequently in the feedback optimization literature \cite[Assumption 2]{bianchin2023online}, \cite[Assumption 1]{colombino2019online}. Under Assumption \ref{ass:cost}(ii), problem \eqref{eq:optimization} admits a unique optimal steady state.

Optimization problems of this form are ubiquitous in many engineering disciplines and are generally challenging to solve using feedforward approaches. This difficulty arises because feedforward methods typically require precise knowledge of the steady-state mapping $h$ and of the disturbance $w$, which is often unavailable in practice. This challenge motivates the development of feedback-based optimization schemes that utilize real-time measurements of the output to eliminate explicit dependence on the disturbance $w$. 

As a starting point for the developments below, we adopt the basic feedback-based projected gradient algorithm, which is the discrete-time counterpart of the continuous-time projected feedback gradient flow proposed in \cite{hauswirth2020timescale}. With a step-size $\alpha > 0$ and projection operator $\text{proj}_{\mathcal{U}}$, the iterations, $k\in\mathbb{N}$, are given by
 \begin{equation}\label{standard_FO}
 \begin{cases}
 y_k = g(x_k,u_k) + w_2\\
 g_k = \nabla_uJ(u_k,y_k) + \nabla h(u_k,w)^{\top}\nabla_yJ(u_k,y_k)\\
 u_{k+1} = \text{proj}_{\mathcal{U}}\left(u_k - \alpha g_k\right)\\
 x_{k+1} = f(x_k,u_k) + w_1.
 \end{cases}\end{equation}
The practical advantage of this approach is that we no longer need a model of the full map $h(u_k,w)$, as this information has been incorporated through $y_k$ in real time. Related feedback-based schemes have already been explored in various problem settings including time-varying optimization problems \cite{colombino2019online}, \cite{bernstein2019online}, distributed implementations for large-scale systems \cite{chang2019saddle}, \cite{carnevale2024nonconvex}, and state-constrained formulations \cite{haberle2020non}, \cite{chen2023online}. However, most existing variants share a common implementation challenge: They require knowledge of the Jacobian matrix $\nabla h(u_k,w)$ with respect to the decision variable $u$, which is difficult to obtain in many real-world applications due to the high dimensionality of the state variables and the complexity of the system model. 

While linearization around a fixed operating point and the linearized input-output mapping offers a tractable approximation, the selection of this linearization point critically affects the performance of the algorithm. As analyzed in \cite{colombino2019towards}, if the linearization point is poorly chosen---namely, far from the true optimal operating point---the resulting approximate algorithm cannot guarantee convergence to a solution close to the globally optimal steady state.

This fundamental limitation motivates the design of the sequential feedback optimization scheme proposed in this paper.

 \section{Sequential FO algorithm and convergence results}\label{section:algorithm}
 \subsection{Algorithm design}
 Instead of fixing the linearization point throughout the optimization process, our proposed sequential approach adaptively updates the linearization along the state and control variables at each time step. By continuously re-linearizing the system around the current operating point as the algorithm progresses, our sequential method can better track the path towards the globally optimal steady state. The key insight is that while a single fixed linearization may be inadequate for global optimization, a sequence of local linearizations that evolve with the optimization trajectory can maintain sufficient accuracy to possibly guide the algorithm towards the true optimum.

We integrate the sequential linearization into the feedback optimization framework by updating the gradient computation at each iteration using the current linearized sensitivity matrix, ensuring local accuracy throughout the optimization trajectory.

Next, we present the sequential FO algorithm.
 \begin{algorithm}
\caption{Sequential Feedback Optimization (SFO)}
\label{alg1}
\textbf{Initialization:} initial state $\hat{x}_0$, initial input $\hat{u}_0\in\mathcal{U}$\\
\textbf{Offline Phase:}
\begin{enumerate} 
\item Design step-size $\alpha$
\end{enumerate}
\textbf{Online Phase (at each time step $k$):}
\begin{enumerate}\setcounter{enumi}{1}
\item Estimate sensitivity: Linearize the dynamical model around $(\hat{x}_{k}, \hat{u}_{k})$ and compute
\begin{align}\label{alg1:sensitivity}
\widetilde{\nabla h}(\hat{x}_{k}, \hat{u}_{k})&= \nabla_x g(\hat{x}_{k}, \hat{u}_{k})(I - \nabla_x f(\hat{x}_{k}, \hat{u}_{k})^{-1}\cdot\notag\\
&\quad\  \nabla_u f(\hat{x}_{k}, \hat{u}_{k}) + \nabla_u g(\hat{x}_{k}, \hat{u}_{k})
\end{align}
\item Output measurement: 
\begin{equation}\label{alg1:output}
\hat{y}_{k} = g\left(\hat{x}_{k},\hat{u}_{k}\right) + w_2
\end{equation}
\item Gradient descent:
\begin{align}\label{alg1:input}
\hat{u}_{k+1} &= \mathrm{proj}_{\mathcal{U}}\Big(\hat{u}_{k} - \alpha\Big(\nabla_uJ\left(\hat{u}_{k}, \hat{y}_{k}\right)\notag\\
&\quad + \widetilde{\nabla h}\left(\hat{x}_{k}, \hat{u}_{k}\right)^{\top}\nabla_yJ\left(\hat{u}_{k}, \hat{y}_{k}\right)\Big)\Big)
\end{align}
\item System dynamics: 
\begin{equation}\label{alg1:state}
\hat{x}_{k+1} = f\left(\hat{x}_{k},\hat{u}_{k}\right) + w_1
\end{equation}
\end{enumerate}
\end{algorithm}
 
 At each time step, the algorithm first linearizes the nonlinear system \eqref{plant_model} around the current operating point $(\hat{x}_{k},\hat{u}_{k})$ to obtain the linearized dynamics
 \begin{equation*}
 \begin{cases}
 \delta x_{k+1} = \nabla_x f(\hat{x}_{k},\hat{u}_{k}) \delta x_k + \nabla_u f(\hat{x}_{k},\hat{u}_{k}) \delta u_k\\
 \delta y_k = \nabla_x g(\hat{x}_{k},\hat{u}_{k}) \delta x_k + \nabla_u g(\hat{x}_{k},\hat{u}_{k}) \delta u_k,
 \end{cases}
 \end{equation*}
 where $\delta x_k = x_k - \hat{x}_{k}$ and $\delta u_k = u_k - \hat{u}_{k}$. The sensitivity estimate $\widetilde{\nabla h}(\hat{x}_{k},\hat{u}_{k})$ in \eqref{alg1:sensitivity} is then computed as the steady-state input-output sensitivity of this linearized system by solving
 \begin{equation*}
 \begin{cases}
 \delta x_{\text{ss}} = (I - \nabla_x f(\hat{x}_{k},\hat{u}_{k}))^{-1}\nabla_u f(\hat{x}_{k},\hat{u}_{k}) \delta u\\
 \delta y_{\text{ss}} = \nabla_x g(\hat{x}_{k},\hat{u}_{k}) \delta x_{\text{ss}} + \nabla_u g(\hat{x}_{k},\hat{u}_{k}) \delta u.
 \end{cases}
 \end{equation*} 
Using this approximate Jacobian matrix, our algorithm updates the control variables through projected gradient descent in combination with the current output measurements, and then applies the updated control to let the system dynamics evolve. Intuitively, as the control variables continue to improve and approach the optimal steady state, Algorithm \ref{alg1} continuously reduces the linearization error through sequential re-linearization, thereby ultimately achieving convergence to an approximation of the optimal steady state.


 \begin{rem}\label{rem:real_time_iteration}
 The sequential linearization approach bears some similarity with real-time iteration schemes in model predictive control (MPC) \cite{diehl2005real}, \cite{gros2020linear}. However, there are fundamental differences: Real-time iteration in MPC aims to track a known reference trajectory or regulate the system to a predetermined steady state, while our sequential FO algorithm seeks to drive the system to an \emph{unknown} optimal steady state defined implicitly by problem \eqref{eq:optimization}. Addtionally, MPC solves finite-horizon dynamic optimization problems, whereas we focus on steady-state optimization using sequential linearizations to approximate the input-output sensitivity matrix.
 \end{rem}

\subsection{Convergence analysis}
In this section, we establish the convergence of Algorithm \ref{alg1}. The analysis proceeds in two stages: we first prove convergence guarantees for the ideal feedback optimization algorithm with exact Jacobian matrix, and then quantify how the approximation errors introduced by sequential linearization affect the convergence behavior.

To focus the analysis on the core dynamics, we assume $g(x,u) = x$, implying that the output directly measures the system state. This setting, which simplifies the notation, is common in the feedback optimization literature \cite{chen2023online}, \cite{hauswirth2020timescale}. The subsequent analysis extends to any output function $g(x, u)$ with a Lipschitz continuous gradient, as the proof structure remains the same, introducing only additional notational complexity.

We begin by characterizing the Lipschitz continuity of the approximate sensitivity $\widetilde{\nabla h}$ defined in \eqref{alg1:sensitivity}, which is instrumental for our subsequent analysis. 
\begin{lem}\label{lem:lipschitz_nablah}
Let Assumption \ref{ass:system} be satisfied. Then the map $(x,u)\mapsto \widetilde{\nabla h}(x,u)$ is Lipschitz continuous over $\mathbb{R}^n\times\mathcal{U}$. Specifically, for all $x_1,x_2 \in\mathbb{R}^n$ and $u_1,u_2 \in\mathbb{R}^p$,
\begin{align}\label{lipschitz_nablah}
&\quad\|\widetilde{\nabla h}(x_1,u_1) - \widetilde{\nabla h}(x_2,u_2)\|\notag\\
& \le \frac{(1 - \rho_f)L_{f,u} + G_u^fL_{f,x}}{(1 - \rho_f)^2}\left(\|x_1 - x_2\| + \|u_1 - u_2\|\right).
\end{align}
\end{lem}

A key feature of FO algorithms is their use of real-time outputs $y_k$ instead of steady-state outputs $h(u_k)$. The resulting error is bounded by the following technical result: 
\begin{lem}\label{lem:output_error}
Under Assumptions \ref{ass:system}-\ref{ass:cost}, define the error $e_k:= \|h(u_k) - y_k\|$. Then it holds that
\begin{equation}\label{lem_proof:e_k}
 e_k \le \rho_f^k\|h(u_0) - y_0\| + \frac{\alpha L_h(G_u^J + L_h G_y^J)}{1 - \rho_f}.
 \end{equation} 
\end{lem}

Next, we establish our convergence result for feedback optimization with perfect gradient information:

\begin{lem}\label{lem:standard_fo}
Let Assumptions \ref{ass:system}-\ref{ass:cost} hold and the step size $\alpha > 0$ be chosen such that the spectral radius $\rho(M)$ of the coefficient matrix
\begin{equation}\label{thm1:coef_matrix}
M:=\begin{bmatrix}\sqrt{1 - 2\alpha\mu_J + \alpha^2L_{J,u}^2} + \alpha C_1 & \alpha C_2\\
 G_u^f & \rho_f\end{bmatrix}
 \end{equation}
 satisfies $\rho(M) < 1$, where
 \begin{equation}\label{thm1:C_constant}
 C_1:= L_hL_{J,y} + G_y^JC,\quad C_2:=L_{J,u} + L_hL_{J,y}
 \end{equation} 
 and $C:= \frac{(1 - \rho_f)L_{f,u} + G_u^fL_{f,x}}{(1 - \rho_f)^2}(1+L_h)$. Then,  the ideal feedback optimization algorithm in \eqref{standard_FO} converges to the optimal steady state $(\bar{u}^{\ast},\bar{y}^{\ast})$ of system \eqref{plant_model}.
\end{lem}

\begin{rem}\label{lem:proof_standard_fo}
Although convergence results for basic feedback optimization schemes have appeared in existing literature, they typically focus on either static plants \cite{colombino2019towards} or continuous-time projected gradient flow for continuous-time dynamic plants \cite{hauswirth2020timescale}. Therefore, the analysis technique may not be directly generalized to the discrete-time version in \eqref{standard_FO}, due to fundamental challenges like step-size restrictions and the introduction of the discretization error.  Lemma \ref{lem:standard_fo} establishes convergence rigorously in a discrete-time setting and the proof technique provides the foundation for subsequent convergence analysis with inexact Jacobian matrix.
\end{rem}

Lemma \ref{lem:standard_fo} provides baseline convergence guarantees when exact gradient information is available. Our sequential algorithm introduces approximation error through linearization, and the key to proving convergence lies in analyzing how this error propagates throughout the optimization dynamics.

Let $(y_k, u_k)$ denote the trajectory of the ideal algorithm with exact gradients in \eqref{standard_FO}, and $(\hat{y}_k, \hat{u}_k)$ denote the trajectory of our sequential FO algorithm. We now state our main convergence result, which focuses on bounding the deviation $\|u_k - \hat{u}_k\|$ between the ideal trajectory and approximated one.

\begin{thm}\label{thm1}
Let Assumptions \ref{ass:system}-\ref{ass:cost} hold. Consider the sequential feedback optimization scheme presented in Algorithm \ref{alg1}, and let the step size $\alpha > 0$ be selected according to the stability condition in Lemma \ref{lem:standard_fo}. Then, the control input sequence $\left(\hat{u}_k\right)_{k\in\mathbb{N}}$ generated by Algorithm \ref{alg1} converges to a neighborhood of the optimal steady state $\bar{u}^*$, with the steady-state error bounded as follows: 
 \begin{equation}
\limsup_{k \to \infty} \|\hat{u}_k - \bar{u}^{\ast}\| \leq \frac{\alpha^2 CL_h(G_u^J + L_h G_y^J)}{(1 - \rho_f)(1 - \rho(M))}.
\end{equation}
\end{thm}

Theorem \ref{thm1} establishes that the use of an approximated sensitivity matrix in Algorithm \ref{alg1} leads to practical convergence, where the control input sequence enters and remains within a neighborhood of the true optimum $\bar{u}^{\ast}$. The size of this neighborhood, as quantified by the asymptotic error bound, is shown to be the order of $\alpha^2$, indicating that a smaller step size generally improves the accuracy of the final solution. In addition, the convergence rate is governed by the spectral radius $\rho(M)$, which depends on both the problem parameters and the step-size choice. A smaller step-size may slow down the convergence rate. Thus, selecting $\alpha$ requires balancing the trade-off between accuracy and convergence speed.

\begin{figure}[h]
	\centering
	\includegraphics[width=.48\textwidth]{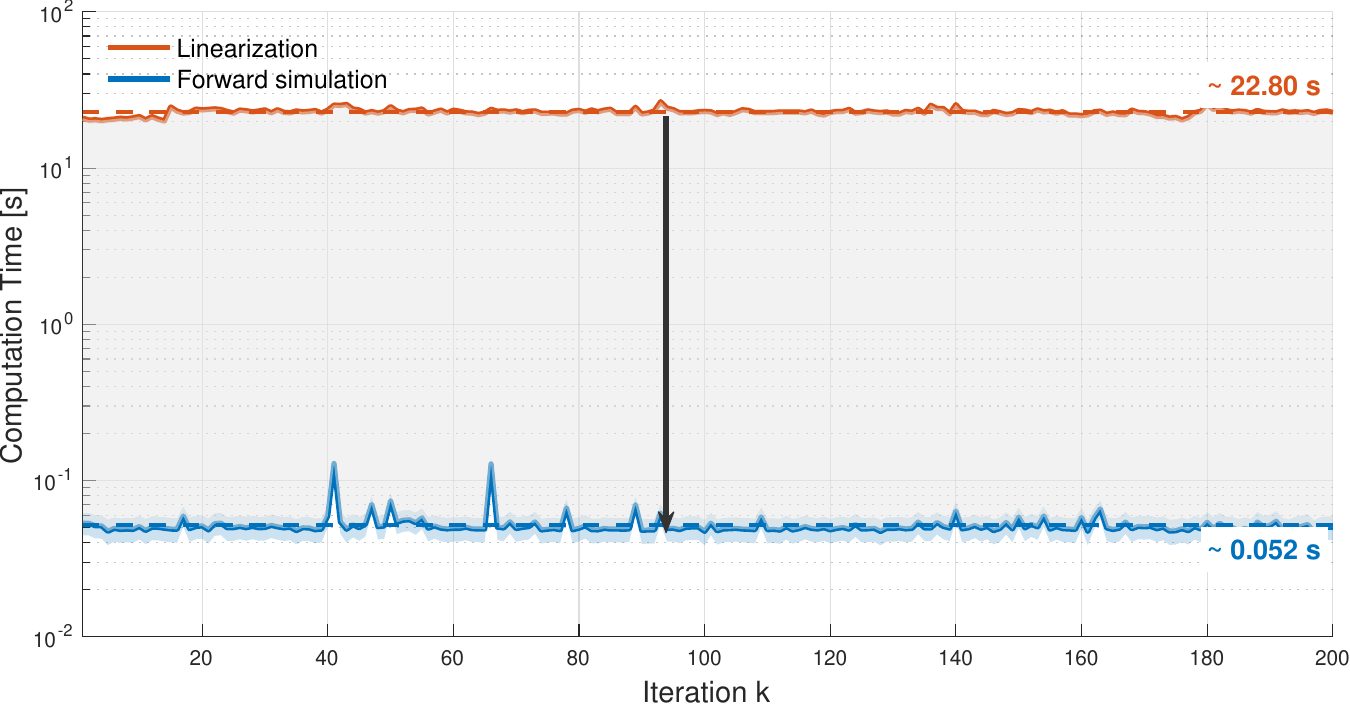}
	\caption{The comparison between linearization time and forward simulation time.}
	\label{FIG:1}
\end{figure}

 \begin{rem}\label{rem:computational_bottleneck}
 For high-dimensional nonlinear dynamic systems, sequential linearization significantly increases computational complexity. In wind farm applications, for instance,  the dynamic flow models typically involve tens of thousands of state variables, making each linearization operation computationally expensive. Empirical studies in Figure \ref{FIG:1} show that for a flow dynamics with around $4.5\times 10^{4}$ states, the linearization time often dominates the forward simulation time, creating a computational bottleneck that motivates a multi-timescale approach.
 \end{rem}

\section{Sequential multi-timescale FO algorithm and convergence results}\label{section:multi-scale}
 To address the computational bottleneck identified above, this section develops a multi-timescale variant that reduces the frequency of linearization operations while maintaining convergence guarantees. The key insight is to separate the linearization updates (slow timescale) from the optimization updates (fast timescale), thereby distributing the computational cost of linearization over multiple optimization steps. Specifically, we further propose a sequential multi-timescale FO algorithm based on Algorithm \ref{alg1}.
 \begin{algorithm}[htbp]
\caption{Sequential Multi-Timescale Feedback Optimization (SMTFO)}
\label{alg:multiscale}
\textbf{Initialization:} initial state $\hat{x}_0$, initial input $\hat{u}_0 \in \mathcal{U}$\\
\textbf{Offline Phase:}
\begin{enumerate}
\item Design step size $\alpha > 0$ and inner loop length $T \geq 1$
\end{enumerate}
\textbf{Online Phase (at each outer iteration $k$):}
\begin{enumerate}\setcounter{enumi}{1}
\item Estimate sensitivity by linearizing around $(\hat{x}_{k},\hat{u}_{k})$:
\begin{align}\label{alg2:sensitivity}
\widetilde{\nabla h}(\hat{x}_{k}, \hat{u}_{k})&= \nabla_x g(\hat{x}_{k}, \hat{u}_{k})(I - \nabla_x f(\hat{x}_{k}, \hat{u}_{k})^{-1}\cdot\notag\\
&\quad\  \nabla_u f(\hat{x}_{k}, \hat{u}_{k}) + \nabla_u g(\hat{x}_{k}, \hat{u}_{k})
\end{align}
\item Inner optimization loop: For $t = 0, 1, \ldots, T-1$:
\begin{enumerate}
\item Output measurement: 
\begin{equation}
\hat{y}_{k,t} = g(\hat{x}_{k,t},\hat{u}_{k,t}) + w_2
\end{equation}
\item Gradient update:
\begin{align}\label{alg2:input}
\hat{u}_{k,t+1} &= \mathrm{proj}_{\mathcal{U}}\Big(\hat{u}_{k,t} - \alpha\Big(\nabla_uJ(\hat{u}_{k,t},\hat{y}_{k,t}) \notag\\
&\quad + \widetilde{\nabla h}(\hat{y}_{k}, \hat{u}_{k})^{\top}\nabla_yJ(\hat{u}_{k,t},\hat{y}_{k,t})\Big)\Big)
\end{align}
\item System update: 
\begin{equation}
\hat{x}_{k,t+1} = f(\hat{x}_{k,t}, \hat{u}_{k,t}) + w_1
\end{equation}
\end{enumerate}
\item Outer loop update:\quad $\hat{u}_{k+1} = \hat{u}_{k,T},\ \hat{x}_{k+1} = \hat{x}_{k,T}$
\end{enumerate}
\end{algorithm}

Algorithm \ref{alg:multiscale} shares structural similarities with sequential quadratic programming (SQP) \cite{nocedal1999numerical}, \cite{byrd2008inexact}, \cite{torrisi2018projected}, where an approximate subproblem is solved at each major iteration. In fact, if we replace the nonlinear output $\hat{y}_{k,t}$ in \eqref{alg2:input} with the linearized output, the inner iteration is equivalent to solving a subproblem, that is, finding the optimal steady state of the linearized system. However, there are important distinctions: In traditional SQP, the goal is to solve the linearized subproblem as accurately as possible (large $T$), whereas in our setting, we must limit $T$ to balance computationally efficiency against the accumulation of linearization errors.

\begin{thm}\label{thm:multiscale_convergence}
Under Assumptions \ref{ass:system}-\ref{ass:cost}, consider Algorithm \ref{alg:multiscale} with parameters $\alpha$ and $T$ such that the spectral radius condition in Lemma \ref{lem:standard_fo} holds. Then the control sequence $\left(\hat{u}_k\right)_{k\in\mathbb{N}}$ generated by Algorithm \ref{alg:multiscale} converges to a neighborhood of $\bar{u}^{\ast}$ with steady-state error bounded as follows:
\begin{align}\label{multiscale:error_bound}
&\quad \limsup_{k \to \infty} \|\hat{u}_k - \bar{u}^{\ast}\|\notag\\
& \leq \frac{\alpha^2 G_y^JC(G_u^J + L_h G_y^J)}{1-\rho(M)}\bigg[\frac{L_h}{1-\rho_f}\notag\\
&\qquad\quad + (T-1)\left(1 + \frac{G_u^f}{1 - \rho_f}\right)\bigg]
\end{align}
where the constant $C$ and the matrix $M$ are defined the same as in Lemma \ref{lem:standard_fo}.
\end{thm}

\begin{rem}
Theorem \ref{thm:multiscale_convergence} quantifies the fundamental trade-off in selecting the inner loop length $T$. A larger $T$ enables faster progress toward the optimal steady state and requiring fewer outer iterations to reach a given neighborhood of the optimum. Since each outer iteration involves costly linearization operations, this faster convergence directly translates to lower overall computational cost. However, it also increases the steady-state error bound by a factor of $(T-1)\left(1 + \frac{G_u^f}{1 - \rho_f}\right)$ due to accumulated linearisation errors. When $T = 1$, the multi-timescale algorithm reduces to Algorithm \ref{alg1}, and the error bound in Theorem \ref{thm:multiscale_convergence} simplifies to that of Theorem \ref{thm1}. This analysis provides guidance for selecting $T$ based on the desired balance between computational cost and solution accuracy in practical implementations.
\end{rem}

\section{Application to wind farm control}\label{section:application}

In this section, we model the wind farm control problem in the feedback optimization framework \eqref{eq:optimization} and use the proposed sequential FO algorithm to solve it. In a wind farm, the wind flowing through upstream wind turbines forms a wake that can significantly lower the extracted power of the downstream wind turbines, thereby decreasing the total wind farm efficiency. Therefore, wind farm control typically aims to maximize the steady-state power production or minimize the fatigue loads by mitigating the wake interaction between the turbines \cite{andersson2021wind}, \cite{vali2019adjoint}.
\subsection{Flow dynamics and problem formulation}
Consider a wind farm composed of $N$ wind turbines, each with rotor radius $R$. For turbine $i$, let $\alpha_i$ denote the axial induction factor and $\gamma_i$ the yaw angle. The yaw angle is defined as the angle between the axial rotor axis and the incoming wind direction. While yaw angle is a control variable that can be directly manipulated, the induction factor is indirectly adjusted by controlling the blade pitch angle and generator torque of the wind turbine \cite{boersma2017tutorial}.

According to the actuator disk model, the power generated by a turbine $i$ is given by \cite[Eq. (5)]{boersma2017tutorial}
\begin{equation}
P_i(\alpha_i, \gamma_i, v_i) = \frac{1}{2}\rho\pi R^2 C_P(\alpha_i, \gamma_i)v_i^3,
\end{equation}
where $C_P(\alpha_i, \gamma_i)$ represents the power coefficient and $v_i$ is the wind speed experienced by the turbine. Due to the wake interaction, $v_i$ is influenced by the control variables of other turbines, i.e., $\alpha_{-i}$ and $\gamma_{-i}$. 

Given this interdependence, maximizing the power production requires individually controlling the turbines in such a way that the collective objective is optimized. Therefore, an understanding of how to model wake interaction is fundamental to developing effective control strategies.

\subsubsection{Medium-fidelity dynamic model (WFSim)}

In this paper, we consider a medium-fidelity wake model that balances computational efficiency with accuracy. Based on the two-dimensional Navier--Stokes (NS) equations, \cite{boersma2018control} developed a dynamic wind flow model, Wind Farm Simulator (WFSim), to effectively capture the wake interaction within a wind farm. This model involves a spatial discretization of the NS equations across a staggered grid, leading to a discrete-time nonlinear dynamics representation. The WFSim can be mathematically expressed as
\begin{equation}\label{WFSim_model}
\begin{cases}
E(X_k)X_{k+1} = AX_k + b(X_k, v_k, \gamma_k),\\
y_{i,k} = P_i(v_{i,k}, \gamma_{i,k}, X_k), \quad i = 1, \ldots, N,
\end{cases}
\end{equation}
where the state variables are defined as $(u_k, v_k, p_k)$, with $u_k$ and $v_k$ representing all the flow velocities in the longitudinal and lateral direction at every point of the staggered grid, and $p_k$ representing the pressure vector. The coefficient matrices are defined by
\begin{equation}
E(X_k) = \begin{bmatrix}
A_u(u_k, v_k) & 0 & B_1 \\
0 & A_v(u_k, v_k) & B_2 \\
B_1^T & 2B_2^T & 0
\end{bmatrix},
\end{equation}
\begin{equation}
A = \begin{bmatrix}
A_{11} & 0 & 0 \\
0 & A_{22} & 0 \\
0 & 0 & 0
\end{bmatrix}.
\end{equation}
The control variables are disk-based thrust coefficients and yaw angles, denoted as
\begin{equation}
v_k = [C_{T_1,k}, C_{T_2,k}, \ldots, C_{T_N,k}]^{\top}
\end{equation}
\begin{equation}
\gamma_k = [\gamma_{1,k}, \gamma_{2,k}, \ldots, \gamma_{N,k}]^{\top}
\end{equation}
where $C_{T_i,k} := \frac{4\alpha_i,k}{1-\alpha_i,k}$ is a function of axial induction factor. This formulation captures the coupling between flow velocities and pressure through the coefficient matrices, while the nonlinear term $b(X_k, v_k, \gamma_k)$ accounts for the actuator disk forces exerted by the turbines. For a detailed description of this dynamical flow model, see \cite{boersma2018control}.

\subsubsection{Wind farm power maximization}

We now formulate the steady-state power maximization problem with the wind flow dynamics. Denoting $y = h(v, \gamma)$ as the steady-state input-output mapping of the dynamic model in \eqref{WFSim_model}, we have
\begin{equation}\label{power_max}
\begin{cases}
\min\limits_{v,\gamma}\  \left(\frac{\mathbf{1}_N^T y - P^{\text{ref}}}{P^{\text{ref}}}\right)^2 + \mu\|v\|^2 + \mu_\gamma\|\gamma\|^2 \\
\text{s.t.} \quad y = h(v, \gamma)\\
\qquad\  v_i \in [v_{\min}, v_{\max}],\ \forall i = 1, \ldots, N \\
\qquad\  \gamma_i \in [\gamma_{\min},\gamma_{\max}],\ \forall i = 1, \ldots, N
\end{cases}
\end{equation}
where $P^{\text{ref}}$ is a given total power reference and $\mu,\mu_\gamma$ are regularization coefficients. The regularization terms are used to prevent excessive control actions and improve numerical stability. Note that the state variable $X_k$ in \eqref{WFSim_model} typically has tens of thousands of elements, thus making it difficult to obtain an explicit form of $h(v, \gamma)$ or $\nabla h(v,\gamma)$. This computational challenge motivates the use of our sequential FO approach.

To overcome the computational challenges and apply Algorithm \ref{alg1}, we derive the linearized sensitivity information. Given an operating point $(X^0, v^0, \gamma^0)$, the linearized WFSim model is obtained by introducing the deviations $\delta X = X - X^0$, $\delta v = v - v^0$, and $\delta \gamma = \gamma - \gamma^0$:
\begin{equation}
E(X^0)\delta X_{k+1} = \mathcal{A}\delta X_k + B_1 \delta v_k + B_2 \delta \gamma_k
\end{equation}
\begin{equation}
\delta y_k = C\delta X_k + D_1 \delta v_k + D_2 \delta \gamma_k
\end{equation}
with the following coefficient matrices:
\begin{equation}
\mathcal{A} = A + \frac{\partial b(X_k, v_k, \gamma_k)}{\partial X_k}\bigg|_{X^0,v^0,\gamma^0} - \frac{\partial E(X_k)X_{k+1}}{\partial X_k}\bigg|_{X^0}
\end{equation}
\begin{equation}
\begin{split}
B_1 &= \frac{\partial b(X_k, v_k, \gamma_k)}{\partial v_k}\bigg|_{X^0,v^0,\gamma^0} \\
B_2 &= \frac{\partial b(X_k, v_k, \gamma_k)}{\partial \gamma_k}\bigg|_{X^0,v^0,\gamma^0}
\end{split}
\end{equation}
\begin{equation}
C = \frac{\partial P(X_k, v_k, \gamma_k)}{\partial X_k}\bigg|_{X^0,v^0,\gamma^0}
\end{equation}
\begin{equation}
\begin{split}
D_1 &= \frac{\partial P(X_k, v_k, \gamma_k)}{\partial v_k}\bigg|_{X^0,v^0,\gamma^0}\\
D_2 &= \frac{\partial P(X_k, v_k, \gamma_k)}{\partial \gamma_k}\bigg|_{X^0,v^0,\gamma^0}.
\end{split}
\end{equation}
Consequently, we get
\begin{equation}
\widetilde{\nabla_v h}(X^0, v^0, \gamma^0) = C(E(X^0) - \mathcal{A})^{-1}B_1 + D_1,
\end{equation}
\begin{equation}
\widetilde{\nabla_\gamma h}(X^0, v^0, \gamma^0) = C(E(X^0) - \mathcal{A})^{-1}B_2 + D_2.
\end{equation}
These linearized input-output sensitivity matrices provide the approximate gradient information required by Algorithm \ref{alg1}, enabling efficient optimization despite the high-dimensional state space. By updating these sensitivities at each iteration, the algorithm adapts to the nonlinear wake dynamics while maintaining computational tractability.

 \subsection{Medium-fidelity simulation results}
 Next, we demonstrate the effectiveness of the proposed sequential feedback optimization strategy through comprehensive simulations using a realistic wind farm configuration. The simulations validate the theoretical convergence results and illustrate the practical benefits of the multi-timescale approach in terms of both power production enhancement and computational efficiency.

\begin{figure}
	\centering
	\includegraphics[width=.48\textwidth]{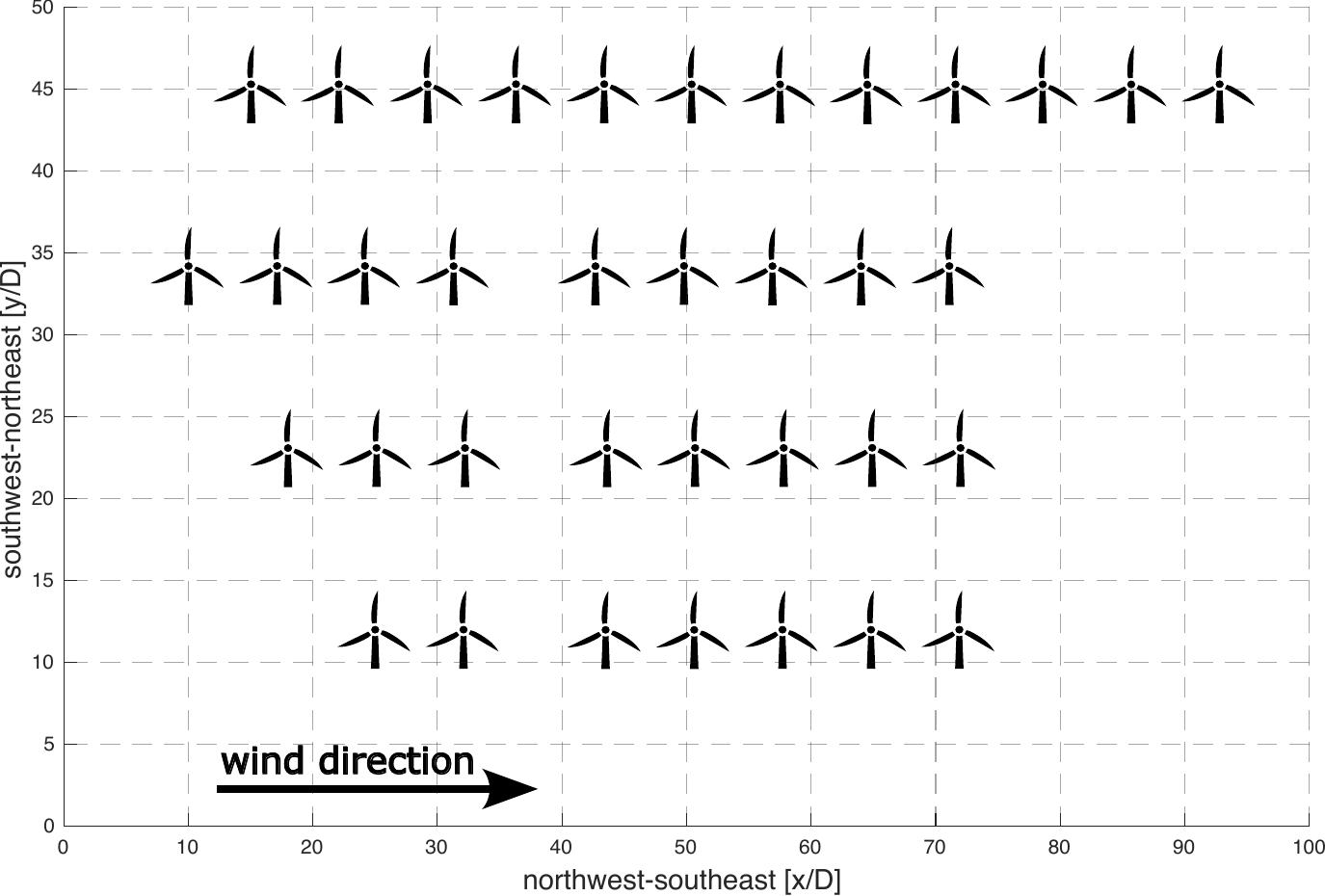}
	\caption{OWEZ wind farm layout.}
	\label{FIG:2}
\end{figure}

The simulations are conducted using the configuration of the Offshore wind farm Egmond aan Zee (OWEZ) \cite{churchfield2015comparison}, \cite{larsen2013validation}, located off the coast of the Netherlands. The farm consists of 36 Vestas V90 3.0 MW wind turbines arranged in a staggered configuration across four rows, as illustrated in Figure \ref{FIG:2}. The inter-row spacing is $11.1$D (where D = 90m is the rotor diameter), while the intra-row turbine spacing is $7.1$D. To accommodate navigation channels, the spacing between turbines 16-17, 24-25, and 31-32 is increased to 11.4D. We have a field of $95\text{D}\times 50\text{D}$ $m^2$ with a staggered grid of $100\times 50$ cells. For the simulations, we assume that the atmospheric conditions are constant and initialize the field with a uniform wind speed of $u = 8~\text{m/s}$ and $v = 0~\text{m/s}$.

The control constraints are set based on physical limitations and operational considerations: The thrust coefficient bounds are $C_{T,\min} = 0.4$ and $C_{T,\max} = 3.6$, corresponding to the feasible operating range of modern wind turbines, while the yaw angles are constrained within $\gamma_{\min} = -30^{\circ}$ and $\gamma_{\max} = 30^{\circ}$ to avoid excessive mechanical stress and maintain grid connectivity. 
 \subsubsection{Comparison with greedy controller}
 The baseline comparison is conducted against the conventional greedy control strategy, which maximizes the power output of each individual turbine without considering wake interactions. In the greedy approach, each turbine operates at its individually optimal thrust coefficient (typically $C_T^{\text{greedy}} = 2$) with zero yaw misalignment ($\gamma^{\text{greedy}} = 0^{\circ}$).

For the sequential FO algorithm, we set the power reference to $P^{\text{ref}} = 36$ MW and choose regularization parameters $\mu = 8\times10^{-4}$ and $\mu_\gamma = 6\times10^{-5}$ to allow sufficient control authority while preventing excessive control actions. The step size is selected as $\alpha = 0.25$ and $\alpha_\gamma = 3$.

\begin{figure}
	\centering
	\includegraphics[width=.48\textwidth]{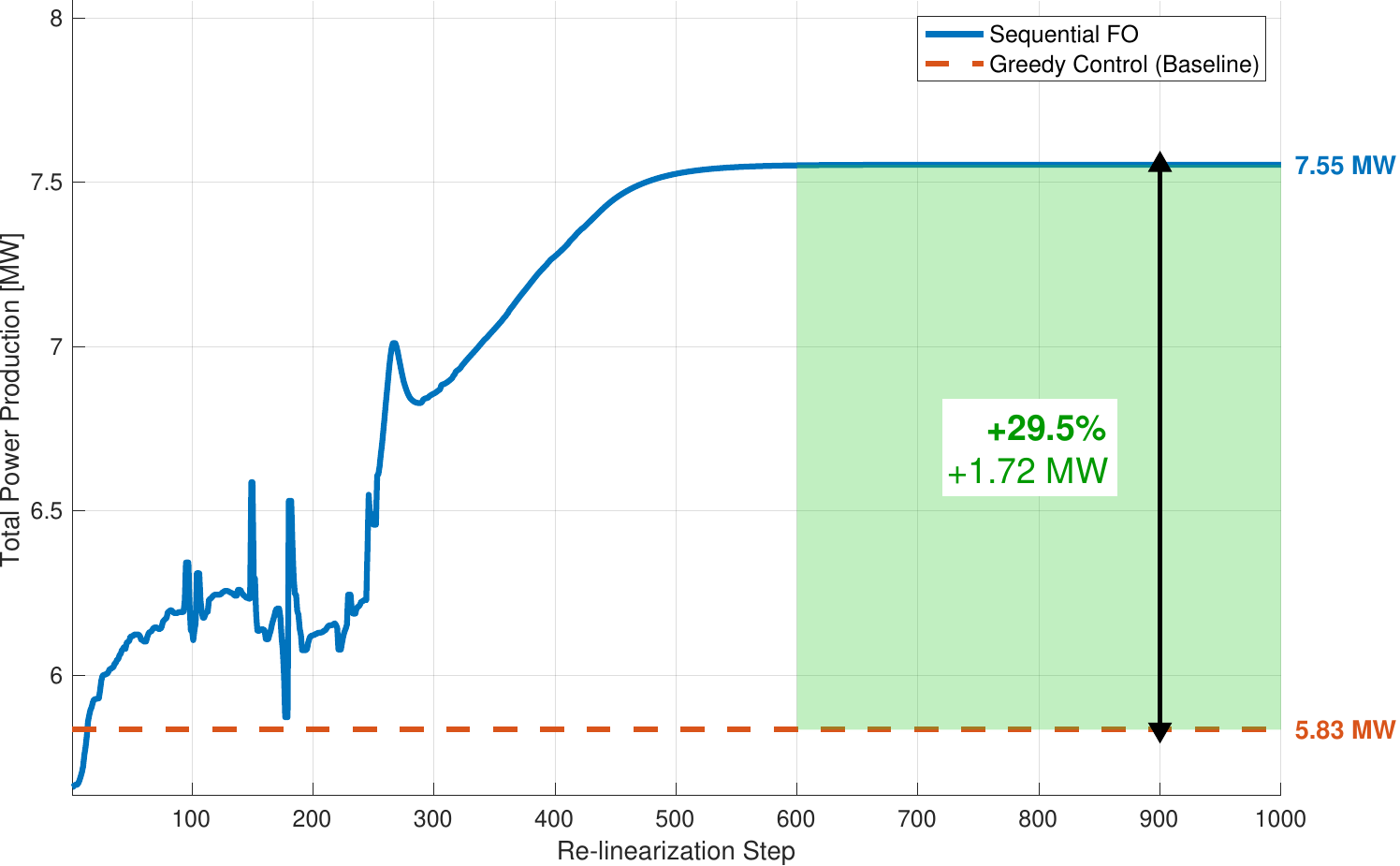}
	\caption{Comparison between sequential FO algorithm and greedy controller.}
	\label{fig:power_comparison}
\end{figure}

Figure \ref{fig:power_comparison} presents the power production trajectories for both control strategies. The results demonstrate that the sequential FO algorithm successfully drives the wind farm to a significantly higher steady-state power output compared to greedy controller. Specifically, the optimized control strategy achieves a $29.46\%$ increase in total power production, rising from approximately $5.83$ MW under greedy control to $7.55$ MW with our sequential FO approach.

This substantial improvement is achieved through coordinated control of upstream turbines, which strategically reduce their individual power extraction to minimize wake losses for downstream turbines. The sequential algorithm automatically discovers this cooperative strategy by leveraging real-time flow measurements and the approximate gradient information obtained through adaptive linearization. 
 \subsubsection{Influence of the inner loop iterations}
 As demonstrated in Remark \ref{rem:computational_bottleneck}, a key advantage of the multi-timescale Algorithm \ref{alg:multiscale} is its ability to reduce computational overhead by performing multiple optimization steps between expensive linearization operations. This subsection investigates the trade-off between computational efficiency and convergence behavior as a function of the inner loop length $T$.

\begin{figure}
	\centering
	\includegraphics[width=.48\textwidth]{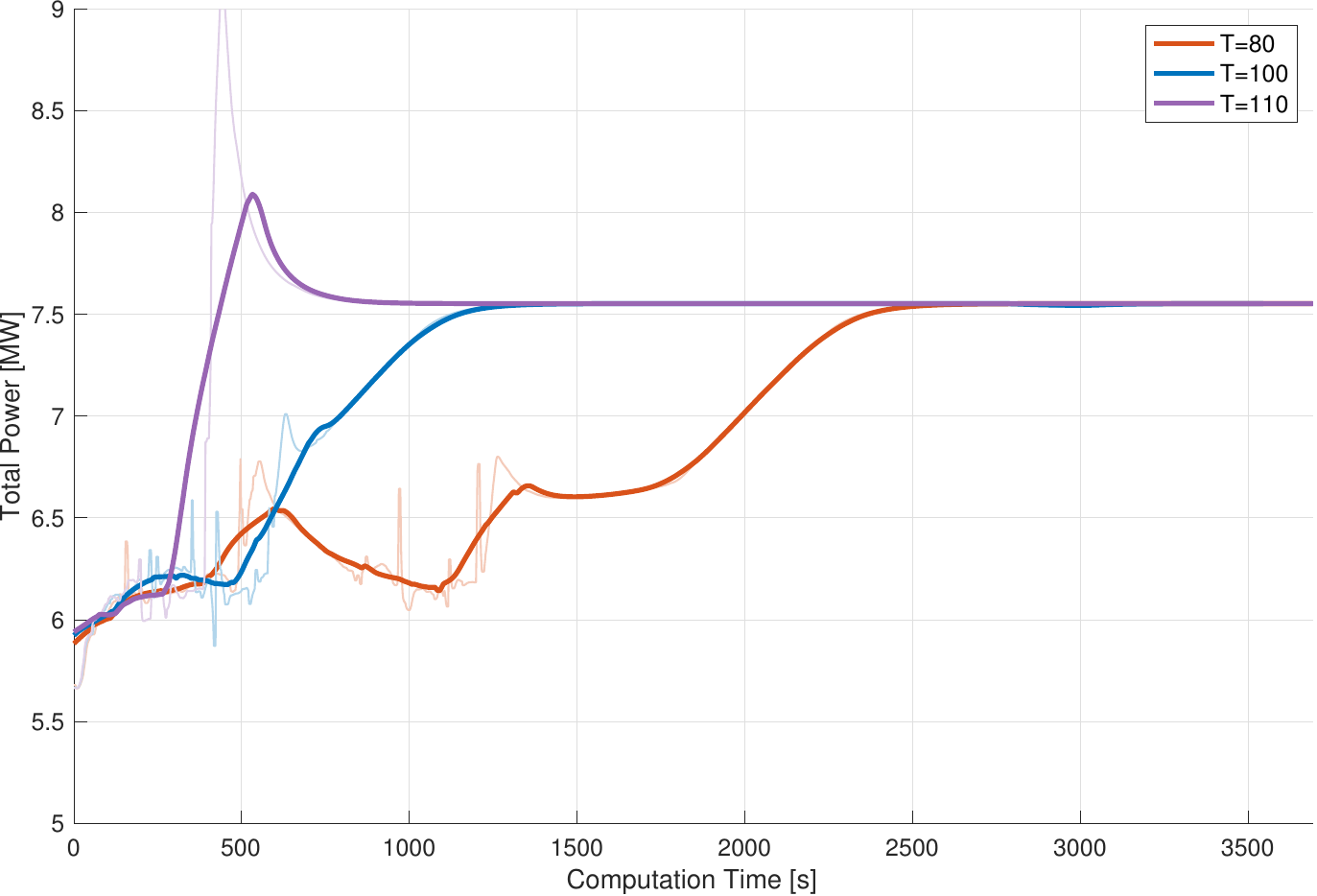}
	\caption{Influence of $T$ on SMTFO. A larger value of $T$ leads to faster convergence vs computation time. For each setting of $T$, the bold line represents the moving average trend of the power trajectory, while the corresponding semi-transparent line shows the unfiltered power trajectory.}
	\label{fig:inner_loop_comparison}
\end{figure}

Figure \ref{fig:inner_loop_comparison} shows the power production trajectories within a given computation time for different values of $T \in \{80, 100, 110\}$. The results confirm the theoretical predictions: Larger values of $T$ enable the algorithm to reach steady state with fewer linearization operations, as each outer iteration allows for more optimization progress.

The numerical results demonstrate the computational benefits of the multi-timescale approach, as detailed in Table \ref{tab:computational_analysis} which provides a detailed breakdown of computational requirements for each configuration, including the number of linearizations, forward simulations, and total computation time on a standard laptop computer. 

\begin{table}[htbp]
\centering
\caption{Computational Time for Different Inner Loop Lengths}\label{tab:computational_analysis}
\setlength{\tabcolsep}{4pt}
\begin{tabularx}{.47\textwidth}{lccc}
\toprule
$T$ & \#Linearizations & \#Forward Sims & Total Time (s)\\
\midrule
80 & 1122 & 89760 & 2400.01  \\
100 & 500 & 50000 & 1178.08   \\
110 & 309 & 33990 & 758.13  \\
\bottomrule
\end{tabularx}
\end{table}

The optimal choice of $T$ represents a balance between computational efficiency and algorithmic stability. For the wind farm configuration studied, $T \in [100,110]$ appears to provide the best compromise, achieving significant computational savings while maintaining robust convergence.
 \subsubsection{Influence of the regularization parameter}
 The regularization parameters $\mu$ and $\mu_\gamma$ in the objective function play a crucial role in balancing power maximization against control effort minimization. This subsection examines their impact on both steady-state performance and control behavior.

Since yaw control typically has a more pronounced effect on power output and wake redirection compared to thrust coefficient adjustments, we focus on the yaw regularization parameter $\mu_\gamma$ while keeping the thrust regularization fixed at $\mu = 8 \times10^{-4}$. The analysis considers values $\mu_\gamma \in \{4\times10^{-5}, 5\times 10^{-5}, 6\times 10^{-5}\}$.

\begin{figure}
	\centering
	\includegraphics[width=.48\textwidth]{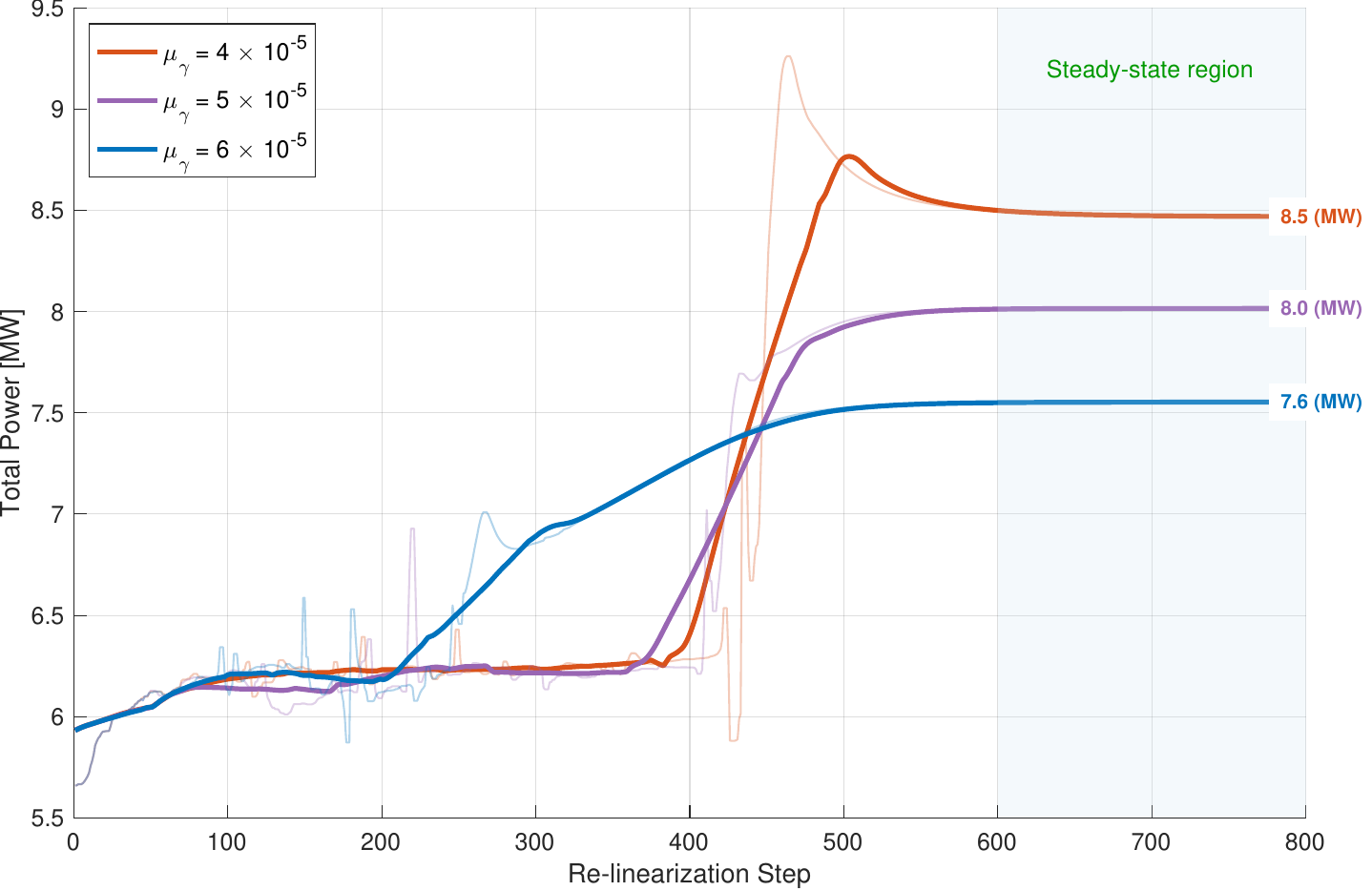}
	\caption{Impact of the yaw regularization parameter $\mu_\gamma$ on the wind farm power production. For each setting of $\mu_\gamma$, the bold line represents the moving average trend of the power trajectory, while the corresponding semi-transparent line shows the unfiltered power trajectory.}
	\label{fig:regularization_power}
\end{figure}
\begin{figure}
	\centering
	\includegraphics[width=.48\textwidth]{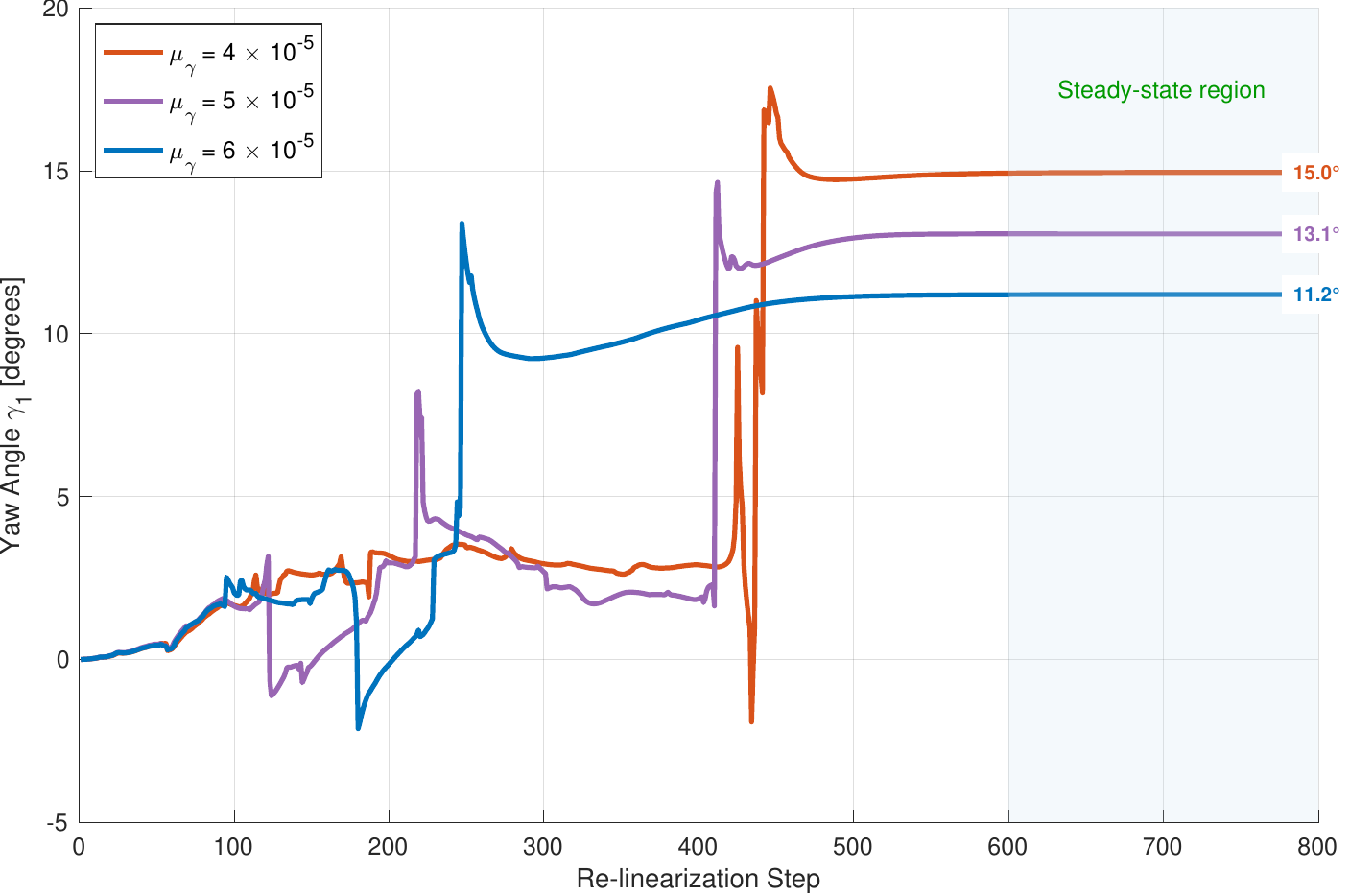}
	\caption{Yaw angle evolution for one of the upstream turbines under different regularization parameters. Smaller $\mu_\gamma$ values result in larger steady-state yaw angles.}
	\label{fig:regularization_yaw}
\end{figure}

Figure \ref{fig:regularization_power} illustrates the relationship between $\mu_\gamma$ and steady-state power production. As expected, reducing $\mu_\gamma$ allows for more aggressive yaw control, leading to higher power output. 
Figure \ref{fig:regularization_yaw} shows the corresponding yaw angles of the most upstream turbine. Smaller values of $\mu_\gamma$ result in larger yaw misalignments for upstream turbines that redirect their wakes to minimize downstream losses. 

While these results demonstrate clear performance benefits with reduced regularization, from a practical perspective, the choice of $\mu_\gamma$ must consider not only power production but also mechanical loading and turbine fatigue. Large yaw angles impose additional structural loads on the nacelle, tower, and foundation systems, potentially reducing operational lifespan and increasing maintenance costs. A formal analysis of the trade-off between power maximization and fatigue load minimization is a promising topic for future research.

\subsection{High-fidelity simulation results}
To further validate the proposed sequential feedback optimization strategy under realistic aerodynamic conditions, we conduct high-fidelity simulations using the Simulator for Offshore Wind Farm Application (SOWFA) \cite{churchfield2012overview}.
\subsubsection{Simulation setup}
Due to computational resource limitations and turbine model availabilty constraints in SOWFA, we consider a simulated wind farm configuration composed of three rows of three NREL 5MW reference wind turbines arranged in a regular grid pattern. The inter-turbine spacing is set to five rotor diameters ($5D$) in the downstream direction and three rotor diameters ($3D$) in the crosswind direction. The accuracy and applicability of the WFSim model for this particular wind farm layout have previously been validated in \cite{boersma2018control}.   

The computational domain in the SOWFA simulation is defined as a $3000~\text{m}\times 1800~\text{m}\times 650~\text{m}$ box discretized with sufficient mesh refinement to accurately capture turbulent wake dynamics. An incoming atmospheric boundary layer is described by a mean wind speed of  $8~\text{m/s}$. In addition, a simulation time step of $0.1~\text{s}$ is employed, and the simulation is run for a total duration of $2000~\text{s}$ to allow sufficient wake development.

In accordance with the algorithm scheme described in Section \ref{section:algorithm}, control actions determined by the WFSim-derived linearized sensitivities are implemented in SOWFA. Subsequently, turbine power outputs obtained from the high-fidelity SOWFA simulations are fed into the algorithm to iteratively update the control inputs.
\subsubsection{Results and analysis}
\begin{figure}
	\centering
	\includegraphics[width=.48\textwidth]{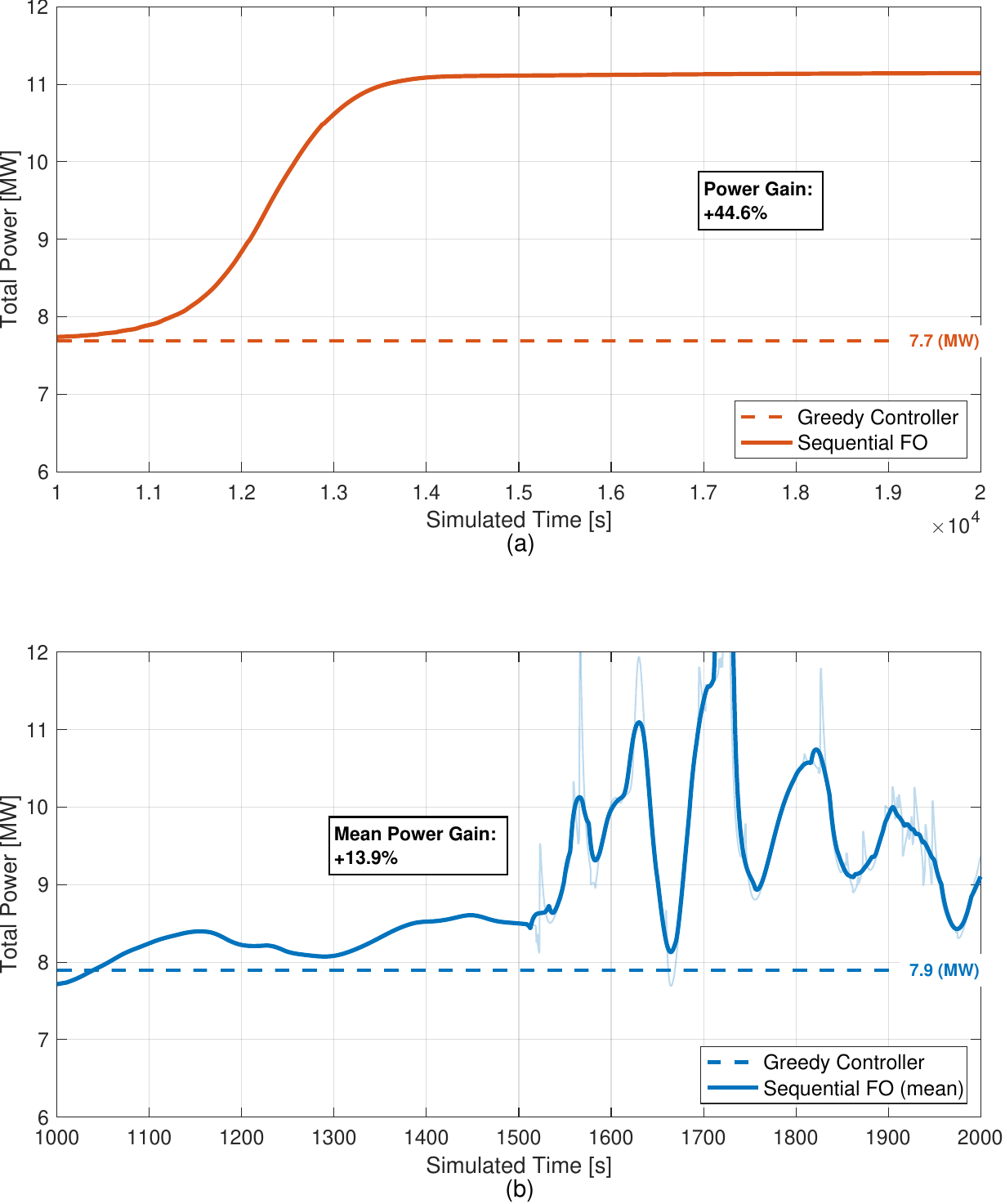}
	\caption{Total power production of the Sequential Feedback Optimization (SFO) strategy. \textbf{(a)} In the nominal WFSim model, SFO demonstrates convergence to a steady state. \textbf{(b)} When validated in the high-fidelity SOWFA model, the same strategy exhibits oscillations due to the inherent model mismatch. The bold line indicates the moving average trend, while the thin semi-transparent line shows the unfiltered power trajectory.}
	\label{fig:high_fidelity}
\end{figure}

The primary validation results are shown in Figure~\ref{fig:high_fidelity}, which compares Algorithm \ref{alg1}'s performance in its nominal design environment (WFSim) with its performance in the high-fidelity SOWFA environment. The initial phase of the SOWFA simulation is treated as a warm-up period to allow the wake interactions to become fully developed. Our analysis focuses on the time window from $1000\text{s}$ to $2000\text{s}$, which represents the algorithm's performance in a realistically waked environment.

As shown in Figure~\ref{fig:high_fidelity}(a), the algorithm exhibits ideal performance in the WFSim environment, achieving convergence to a steady state with an approximate power gain of $44.58\%$. In contrast, Figure~\ref{fig:high_fidelity}(b) reveals the challenge of applying this model-based strategy in a more realistic environment. The power oscillations are a direct consequence of model mismatch. WFSim employs a two-dimensional representation of the Navier--Stokes equations, which simplifies vertical flow structures, turbulent mixing phenomena, and wake recovery processes critical to accurately modelling turbine interactions. SOWFA, on the other hand, incorporates a full three-dimensional Large Eddy Simulation (LES) approach, capturing detailed turbulent structures and more realistic wake dynamics. In fact, such performance gap when transferring controllers from lower- to high-fidelity models has also been observed and discussed in prior wind farm control literature \cite{santoni2015development}, \cite{gebraad2016wind}, \cite{annoni2016analysis}. 

Despite the lack of stable convergence, the result from SOWFA is still valuable. By calculating the moving average of the power trajectory (bold blue line in Figure~\ref{fig:high_fidelity}(b)), we observe a mean power gain of $13.88\%$ over the greedy baseline. It demonstrates that the proposed sequential feedback optimization algorithm effectively steers the wind farm to a new operating region with higher average power output. This validates the practical potential of the feedback optimization framework for improving wind farm operational performance under realistic conditions.

\section{Conclusion}\label{section:conclusion}
This paper demonstrates that adaptive linearization within feedback optimization effectively addresses the fundamental challenge of computing accurate sensitivity information for nonlinear dynamical systems. By sequentially updating the linearization point as the system evolves, our approach avoids the dependence on a fixed operating point, thus improving the optimization accuracy. The developed multi-timescale implementation further reduces computational complexity while maintaining theoretical convergence guarantees. Specifically, our analysis explicitly quantifies the relationship between computational efficiency and steady-state accuracy. Numerical results from realistic wind farm control simulations confirm the theoretical findings and show significantly improved performance compared to the greedy control strategy, indicating the potential applicability of our proposed methods for engineering systems.

Future research will explore extensions of this framework to nonlinear dynamics influenced by stochastic disturbances and the development of distributed implementations suitable for large-scale networked systems.
\begin{appendix}
\section*{Appendix A. Proof of Lemma~\ref{lem:lipschitz_nablah}}
\renewcommand{\theequation}{A.\arabic{equation}}
\setcounter{equation}{0}
Since $g(x,u) = x$, the definition in \eqref{alg1:sensitivity} simplifies to
\begin{equation}\label{simplified_sensitivity}
\widetilde{\nabla h}(x_1,u_1) = (I - \nabla_x f(x_1,u_1)^{-1}\nabla_u f(x_1,u_1).
\end{equation}
Therefore, we derive
\begin{align}\label{lem_proof:gradient_h_error}
 &\quad \left\|\widetilde{\nabla h}(x_1,u_1) - \widetilde{\nabla h}(x_2,u_2)\right\|\notag\\
  &= \Big\|(I - \nabla_x f(x_1,u_1))^{-1}\nabla_u f(x_1,u_1)\notag\\
  &\quad\ \  - (I - \nabla_x f(x_2,u_2)^{-1}\nabla_u f(x_2,u_2)\Big\|\notag\\
 &= \Big\|(I - \nabla_x f(x_1,u_1))^{-1}[\nabla_u f(x_1,u_1) - \nabla_u f(x_2,u_2)]\notag\\
 &\quad\ \  + (I - \nabla_x f(x_1,u_1))^{-1}\nabla_u f(x_2,u_2)\notag\\
 &\quad\ \  - (I - \nabla_x f(x_1,u_1))^{-1}\nabla_u f(x_2,u_2)\Big\|\notag\\
 &\le L_{f,u}\left\|(I - \nabla_x f(x_1,u_1))^{-1}\right\|\left(\|x_1 - x_2\| + \|u_1 - u_2\|\right)\notag\\
 &\quad + G_u^f\left\|(I - \nabla_x f(x_1,u_1))^{-1} - (I - \nabla_x f(x_2,u_2))^{-1}\right\|,
 \end{align}
 where the last inequality follows by Assumption \ref{ass:system}(ii). Furthermore, by Assumption \ref{ass:system}(i),
 \begin{equation}\label{matrix_relation1}
 \left\|(I - \nabla_x f(x_1,u_1))^{-1}\right\| \le \sum_{i=0}^{\infty} \|\nabla_x f(x_1,u_1)\|^i \le \frac{1}{1 - \rho_f}.
 \end{equation}
 Note that we have the relation $A^{-1} - B^{-1} = A^{-1}(B-A)B^{-1}$ for any invertible matrices $A$ and $B$. Thus,
 \begin{align}\label{matrix_relation2}
 &\quad\left\|(I - \nabla_x f(x_1,u_1))^{-1} - (I - \nabla_x f(x_2,u_2))^{-1}\right\|\notag\\
 &= \Big\|(I - \nabla_x f(x_1,u_1))^{-1} (\nabla_x f(x_1,u_1) - \nabla_x f(x_2,u_2))\cdot\notag\\
 &\quad\quad (I - \nabla_x f(x_2,u_2))^{-1}\Big\|\notag\\
 &\le \frac{1}{(1-\rho_f)^2}\|\nabla_x f(x_1,u_1) - \nabla_x f(x_2,u_2)\|\notag\\
 &\le \frac{L_{f,x}}{(1-\rho_f)^2}\left(\|x_1 - x_2\| + \|u_1 - u_2\|\right)
 \end{align}
 By substituting \eqref{matrix_relation1} and \eqref{matrix_relation2} into \eqref{lem_proof:gradient_h_error}, we obtain \eqref{lipschitz_nablah}. 
 \section*{Appendix B. Proof of Lemma~\ref{lem:output_error}}
\renewcommand{\theequation}{B.\arabic{equation}}
\setcounter{equation}{0}
By Assumption \ref{ass:system}(iii), we have
\begin{align*}
 e_{k+1} &= \|h(u_{k+1}) - y_{k+1}\|\\
 &\le \|h(u_{k+1}) - h(u_k)\| + \|h(u_{k}) - y_{k+1}\|\\
 &\le L_h\|u_{k+1} - u_k\| + \|f(h(u_k), u_k) - f(y_k, u_k)\|\\
 &\le L_h\|u_{k+1} - u_k\| + \rho_f e_k.
 \end{align*}
 Consequently, via this recursive relationship, we derive
 \begin{align}\label{lem_proof:e_k_imme}
 e_k&\le \rho_f^k\|h(u_0) - y_0\| + \sum_{j=0}^{k-1}\rho_f^{k-1-j}L_h\|u_{j+1} - u_j\|.
 \end{align}
 In addition, from the update rule in \eqref{standard_FO}, we have
 \begin{align}\label{u_consec}
 \|u_{j+1} - u_j\| &\le \alpha \left\|\nabla_uJ(u_j,y_j) + \nabla h(u_j)^{\top}\nabla_yJ(u_j,y_j)\right\|\notag\\
 &\le \alpha (G_u^J + L_h G_y^J).
 \end{align}
 Combining \eqref{lem_proof:e_k_imme} and \eqref{u_consec} yields \eqref{lem_proof:e_k}.
 
 \section*{Appendix C. Proof of Lemma~\ref{lem:standard_fo}}
\renewcommand{\theequation}{C.\arabic{equation}}
\setcounter{equation}{0}
Since $(\bar{u}^{\ast}, \bar{y}^{\ast})$ is the optimal solution of the problem in \eqref{eq:optimization}, the optimality condition gives
\[\inf_{u\in\mathcal{U}}\nabla\tilde{J}(\bar{u}^{\ast})^{\top}(u - \bar{u}^{\ast}) \ge 0,\]
where $\nabla\tilde{J}(u)$ is the gradient of the reduced-form objective function $\tilde{J}(u) := J(u, h(u,w))$. This condition is equivalent to the fixed-point equation in \eqref{standard_FO} for any $\alpha > 0$:
\begin{equation}\label{opt_condition}
\bar{u}^{\ast} \in \text{proj}_{\mathcal{U}}\left(\bar{u}^{\ast} - \alpha \nabla_u\tilde{J}(\bar{u}^{\ast})\right).
\end{equation}
By the non-expansiveness property of the projection operator and by using \eqref{opt_condition}, we get 
 \begin{align}\label{lem_proof:u_error}
 &\quad \|u_{k+1} - \bar{u}^{\ast}\|\notag\\
 &\le\Big\|u_k - \bar{u}^{\ast} - \alpha(\nabla_uJ(u_k,y_k) - \nabla_uJ(\bar{u}^{\ast},\bar{y}^{\ast}))\notag\\
 &\quad\ \  - \alpha(\nabla h(u_k)^{\top}\nabla_yJ(u_k,y_k) - \nabla h(\bar{u}^{\ast})^{\top}\nabla_yJ(\bar{u}^{\ast},\bar{y}^{\ast}))\Big\|\notag\\
 &= \Big\|u_k - \bar{u}^{\ast} - \alpha\left(\nabla_uJ(u_k,y_k) - \nabla_uJ(\bar{u}^{\ast}, y_k)\right)\notag\\
 &\quad\ \  -\alpha\left(\nabla_uJ(\bar{u}^{\ast}, y_k) - \nabla_uJ(\bar{u}^{\ast},\bar{y}^{\ast})\right)\notag\\
 &\quad\ \  - \alpha(\nabla h(u_k)^{\top}\nabla_yJ(u_k,y_k) - \nabla h(\bar{u}^{\ast})^{\top}\nabla_yJ(\bar{u}^{\ast},\bar{y}^{\ast}))\Big\|\notag\\
 &\le \underbrace{\left\|u_k - \bar{u}^{\ast} - \alpha\left(\nabla_uJ(u_k,y_k) - \nabla_uJ(\bar{u}^{\ast}, y_k)\right)\right\|}_{\text{Term 1}}\notag\\
 &\quad + \alpha\left\|\nabla_uJ(\bar{u}^{\ast}, y_k) - \nabla_uJ(\bar{u}^{\ast},\bar{y}^{\ast})\right\|\notag\\
 &\quad + \alpha\underbrace{\left\|\nabla h(u_k)^{\top}\nabla_yJ(u_k,y_k) - \nabla h(\bar{u}^{\ast})^{\top}\nabla_yJ(\bar{u}^{\ast},\bar{y}^{\ast})\right\|}_{\text{Term 2}}.
 \end{align} 
 Term 1 can be bounded by leveraging the strong monotonicity and Lipschitz continuity of $J(u, y)$ with respect to its first argument:
 \begin{align}\label{lem_proof:term1}
 (\text{Term 1})^2 &= \left(u_k - \bar{u}^{\ast} - \alpha(\nabla_uJ(u_k,y_k) - \nabla_uJ(\bar{u}^{\ast},y_k))\right)^{\top}\notag\\
 &\quad \left(u_k - \bar{u}^{\ast} - \alpha(\nabla_uJ(u_k,y_k) - \nabla_uJ(\bar{u}^{\ast},y_k))\right)\notag\\
 &= \|u_k - \bar{u}^{\ast}\|^2 + \alpha^2\left\|\nabla_uJ(u_k,y_k) - \nabla_uJ(\bar{u}^{\ast},y_k))\right\|^2\notag\\
 &\quad - 2\alpha(\nabla_uJ(u_k,y_k) - \nabla_uJ(\bar{u}^{\ast},y_k))^{\top}(u_k - \bar{u}^{\ast})\notag\\
 &\le (1 + \alpha^2L_{J,u}^2 - 2\alpha\mu_J)\|u_k - \bar{u}^{\ast}\|^2,
 \end{align}
 Moreover, for Term 2, applying Assumption \ref{ass:cost}(i) yields
 \begin{align}\label{lem_proof:term2}
 \text{Term 2} &= \Big\|\nabla h(u_k)^{\top}(\nabla_yJ(u_k,y_k) - \nabla_yJ(\bar{u}^{\ast},\bar{y}^{\ast}))\notag\\
 &\quad + (\nabla h(u_k) - \nabla h(\bar{u}^{\ast}))^{\top}\nabla_yJ(\bar{u}^{\ast},\bar{y}^{\ast})\Big\|\notag\\
 &\le L_h L_{J,y}\left(\|u_k - \bar{u}^{\ast}\| + \|y_k - \bar{y}^{\ast}\|\right)\notag\\
 &\quad + G_y^J\left\|\nabla h(u_k) - \nabla h(\bar{u}^{\ast})\right\|.
 \end{align}
 To bound the term $\|\nabla h(u_k) - \nabla h(\bar{u}^{\ast})\|$, we note that the steady-state map is defined implicitly by the fixed-point relation
 \[h(u_k) = f(h(u_k), u_k) + w_1,\]
 which implies
 \[\nabla h(u_k) = \nabla_u f(h(u_k), u_k) + \nabla_x f(h(u_k),u_k)\nabla h(u_k).\]
 By rearranging the terms to solve for $\nabla h(u_k)$ and recalling the definition of the estimated sensitivity matrix $\widetilde{\nabla h}$ in \eqref{simplified_sensitivity}, we obtain
 \begin{align}\label{tilde_h}
 \nabla h(u_k) &= (I - \nabla_x f(h(u_k), u_k))^{-1}\nabla_u f(h(u_k), u_k)\notag\\
 &= \widetilde{\nabla h} (h(u_k), u_k).
 \end{align}
 Next, by Lemma \ref{lem:lipschitz_nablah} and Assumption \ref{ass:system}(iii), we have
 \begin{align}\label{lem1_proof:gradient_h_error}
 &\quad \left\|\nabla h(u_k) - \nabla h(\bar{u}^{\ast})\right\|\notag\\
 &= \frac{(1 - \rho_f)L_{f,u} + G_u^fL_{f,x}}{(1 - \rho_f)^2}\left(\|u_k - \bar{u}^{\ast}\| + \|h(u_k) - h(\bar{u}^{\ast})\|\right)\notag\\
 &\le \frac{(1 - \rho_f)L_{f,u} + G_u^fL_{f,x}}{(1 - \rho_f)^2} (1+ L_h) \|u_k - \bar{u}^{\ast}\|.
 \end{align}
 By substituting \eqref{lem_proof:term1}-\eqref{lem1_proof:gradient_h_error} back into \eqref{lem_proof:u_error} and collecting terms, we derive 
 \begin{align}\label{lem_proof:u_final_error}
 \|u_{k+1} - \bar{u}^{\ast}\| &\le (\sqrt{1 - 2\alpha\mu_J + \alpha^2L_{J,u}^2} + \alpha C_1)\|u_k - \bar{u}^{\ast}\|\notag\\
 &\quad + \alpha C_2\|y_k - \bar{y}^{\ast}\|.\end{align}
 with constants $C, C_1, C_2$ defined in \eqref{thm1:C_constant}.
 
 Next, we analyze the evolution of the error $\|y_{k+1} - \bar{y}^{\ast}\|$. Using the system dynamics $y_{k+1} = f(y_k,u_k) + w_1$ and the steady-state condition $\bar{y}^{\ast} = f(\bar{y}^{\ast},\bar{u}^{\ast}) + w_1$, and applying Assumption \ref{ass:system}, we obtain:
 \begin{align}\label{lem_proof:y_final_error}
 \|y_{k+1} - \bar{y}^{\ast}\| &= \|f(y_k, u_k) - f(\bar{y}^{\ast}, \bar{u}^{\ast})\|\notag\\
 &\le G_u^f\|u_k - \bar{u}^{\ast}\| + \rho_f\|y_k - \bar{y}^{\ast}\|.\end{align}
 Finally, we combine the derived bounds in \eqref{lem_proof:u_final_error} and \eqref{lem_proof:y_final_error} to get the linear matrix inequality
 \[\begin{bmatrix}\|u_{k+1} - \bar{u}^{\ast}\|\\
 \|y_{k+1} - \bar{y}^{\ast}\|\end{bmatrix}\le M\begin{bmatrix}\|u_{k} - \bar{u}^{\ast}\|\\
 \|y_{k} - \bar{y}^{\ast}\|\end{bmatrix} \]
 with
 \[M:=\begin{bmatrix}\sqrt{1 - 2\alpha\mu_J + \alpha^2L_{J,u}^2} + \alpha C_1 & \alpha C_2\\
 G_u^f & \rho_f\end{bmatrix}.\]
 Since $\alpha$ is chosen such that $\rho(M)<1$, the error vector converges linearly to zero, which implies that $\lim_{k\to\infty}\|u_k - \bar{u}^{\ast}\| = 0$ and $\lim_{k\to\infty}\|y_k - \bar{y}^{\ast}\| = 0$. This completes the proof of convergence to the optimal steady state.
 
 \section*{Appendix D. Proof of Theorem~\ref{thm1}}
\renewcommand{\theequation}{D.\arabic{equation}}
\setcounter{equation}{0}
Following a similar analytical structure as in the proof of Lemma \ref{lem:standard_fo}, we first bound the evolution of the input error $\|u_k - \hat{u}_k\|$. The key difference arises from the use of the approximated sensitivity $\widetilde{\nabla h}(\hat{y}_{k}, \hat{u}_{k})$, which introduces an additional term into the analysis: 
  \begin{align}\label{proof:u_error}
 &\quad \|u_{k+1} - \hat{u}_{k+1}\|\notag\\
  &\le (\sqrt{1 - 2\alpha\mu_J + \alpha^2L_{J,u}^2} + \alpha L_hL_{J,y})\|u_k - \hat{u}_k\|\notag\\
&\quad + \alpha (L_{J,u} + L_hL_{J,y})\|y_k - \hat{y}_k\|\notag\\
&\quad + \alpha G_y^J\left\|\nabla h(u_k) - \widetilde{\nabla h}(\hat{y}_{k}, \hat{u}_{k})\right\|.\end{align}
The final term in \eqref{proof:u_error} quantifies the discrepancy in the sensitivity approximation. From \eqref{tilde_h}, we apply Lemma \ref{lem:lipschitz_nablah} to bound this discrepancy as follows:
\begin{align}\label{proof:gradient_h_error}
 &\quad \left\|\nabla h(u_k) - \widetilde{\nabla h}(\hat{x}_{k}, \hat{u}_{k})\right\|\notag\\
 &= \frac{(1 - \rho_f)L_{f,u} + G_u^fL_{f,x}}{(1 - \rho_f)^2}\left(\|u_k - \hat{u}_k\| + \|h(u_k) - \hat{y}_k\|\right).
 \end{align}
 Moreover, we note that
 \begin{align}\label{proof:h_error}
 \|h(u_k) - \hat{y}_k\| &\le \|h(u_k) - h(\hat{u}_k)\| + \|\hat{y}_k - h(\hat{u}_k)\|\notag\\
 &\le L_h \|u_k - \hat{u}_k\| + \|\hat{y}_k - h(\hat{u}_k)\|.
 \end{align}
 By defining $\hat{e}_k := \|\hat{y}_k - h(\hat{u}_k)\|$ as the instantaneous output-measurement mismatch and using analysis analogous to that in Lemma \ref{lem:output_error}, we derive
 \begin{equation}\label{proof:e_k}
 \hat{e}_k \le \rho_f^k\|h(u_0) - y_0\| + \frac{\alpha L_h(G_u^J + L_h G_y^J)}{1 - \rho_f}.
 \end{equation} 
 Substituting \eqref{proof:gradient_h_error} and \eqref{proof:h_error} back into \eqref{proof:u_error} yields a final bound on the input error dynamics: 
 \begin{align}\label{proof:u_final_error}
 &\quad \|u_{k+1} - \hat{u}_{k+1}\|\notag\\
  &\le (\sqrt{1 - 2\alpha\mu_J + \alpha^2L_{J,u}^2} + \alpha C_1)\|u_k - \hat{u}_k\|\notag\\
 &\quad + \alpha C_2\|y_k - \hat{y}_k\| + \alpha G_y^JC\hat{e}_k\end{align}
 with constants $C, C_1, C_2$ defined in \eqref{thm1:C_constant}.
 
 The dynamics of the output error $\|y_k - \hat{y}_k\|$ remain unchanged from the previous analysis:
 \begin{align}\label{proof:y_final_error}
 \|y_{k+1} - \hat{y}_{k+1}\| &= \|f(y_k, u_k) - f(\hat{y}_k, \hat{u}_k)\|\notag\\
 &\le G_u^f\|u_k - \hat{u}_k\| + \rho_f\|y_k - \hat{y}_k\|.\end{align}
 By combining \eqref{proof:u_final_error} and \eqref{proof:y_final_error}, we obtain a linear matrix inequality for the error dynamics:
 \[\begin{bmatrix}\|u_{k+1} - \hat{u}_{k+1}\|\\
 \|y_{k+1} - \hat{y}_{k+1}\|\end{bmatrix}\le M\begin{bmatrix}\|u_{k} - \hat{u}_{k}\|\\
 \|y_{k} - \hat{y}_{k}\|\end{bmatrix} + \alpha G_y^JC\begin{bmatrix}\hat{e}_k\\
 0\end{bmatrix}, \]
 where $M$ is the contraction matrix defined in \eqref{thm1:coef_matrix}. Since $\rho(M)<1$, the homogeneous system is stable. The persistent forcing term, which is asymptotically bounded by \eqref{proof:e_k}, ensures that the error converges to a neighborhood of the origin. Specifically, we obtain
 \begin{align*}
 \limsup_{k \to \infty}\|u_k - \hat{u}_k\|&\le \frac{\alpha C}{1-\rho(M)}\limsup_{k \to \infty}e_k\\
 &\le \frac{\alpha^2 CL_h(G_u^J + L_h G_y^J)}{(1 - \rho_f)(1 - \rho(M))}.
 \end{align*}
 The proof then follows by applying Lemma \ref{lem:standard_fo}.
 
 \section*{Appendix E. Proof of Theorem~\ref{thm:multiscale_convergence}}
\renewcommand{\theequation}{E.\arabic{equation}}
\setcounter{equation}{0}
For each $k$, let us analyze the error $\|\hat{u}_{k, t+1} - u_{kT + t + 1}\|$. With an analysis similar to \eqref{proof:u_error}, we get 
 \begin{align}\label{multiscale:u_error}
 &\quad \|\hat{u}_{k, t+1} - u_{kT + t + 1}\|\notag\\
 &\le \sqrt{1 + \alpha^2L_{J,u}^2 - 2\alpha\mu_J}\|\hat{u}_{k,t} - u_{kT + t}\|\notag\\
 &\quad + \alpha L_{J,u}\|\hat{y}_{k,t} - y_{kT + t}\|\notag\\
 &\quad + \alpha L_hL_{J,y}\left(\|\hat{u}_{k,t} - u_{kT + t}\| + \|\hat{y}_{k,t} - y_{kT + t}\|\right)\notag\\
 &\quad + \alpha G_y^J\left\|\nabla h(u_{kT + t}) - \widetilde{\nabla h}\left(\hat{x}_{k},\hat{u}_{k}\right)\right\|.
 \end{align}
 For the last term, by exploiting the Lipschitz continuity of $\widetilde{\nabla h}$, we derive
 \begin{align}\label{multiscale:gradient_h_error}
 &\quad\|\nabla h(u_{kT + t}) - \widetilde{\nabla h}(\hat{x}_{k}, \hat{u}_{k})\|\notag\\ 
 &\le \|\nabla h(u_{kT + t}) - \widetilde{\nabla h}(\hat{x}_{k,t}, \hat{u}_{k,t})\|\notag\\
 &\quad + \|\widetilde{\nabla h}(\hat{x}_{k,t}, \hat{u}_{k,t}) - \widetilde{\nabla h}(\hat{y}_{k}, \hat{u}_{k})\|\notag\\
 &\le C(\|u_{kT + t} - \hat{u}_{k,t}\| + \|h(u_{kT + t}) - \hat{y}_{k,t}\|)\notag\\
 &\quad + C(\|\hat{u}_{k,t} - \hat{u}_k\| + \|\hat{y}_{k,t} - \hat{y}_k\|).
  \end{align}
 From the update rule in \eqref{alg2:input}, $\forall t \in \{0,1,\dots, T-1\}$: 
 \begin{equation}\label{multiscale:u_cum_error}
 \|\hat{u}_{k,t} - \hat{u}_k\| \le \alpha (T-1)(G_u^J + L_h G_y^J).
 \end{equation}
 Furthermore, we define $\hat{e}_{kT + t} := \|\hat{y}_{k,t+1} - \hat{y}_{k,t}\|$ and derive
 \begin{align*}
 \|\hat{e}_{kT + t}\| &= \|f(\hat{y}_{k,t}, \hat{u}_{k,t}) - f(\hat{y}_{k,t-1}, \hat{u}_{k,t-1})\|\\
 &\le \rho_f \|\hat{y}_{k,t} - \hat{y}_{k,t-1}\| + G_u^f\|\hat{u}_{k,t} - \hat{u}_{k,t-1}\|\\
 &\le \rho_f \|\hat{e}_{kT + t - 1}\| + \alpha G_u^f(G_u^J + L_h G_y^J).
 \end{align*}
 Unrolling the recursion, we get
 \begin{equation*}
 \|\hat{e}_{kT + t}\| \le \rho_f^{kT + t} \|\hat{e}_{0}\| + \frac{\alpha G_u^f(G_u^J + L_h G_y^J)}{1-\rho_f}.
 \end{equation*}
  Therefore, for all $t = 0,1,\dots, T-1$, we have
  \begin{equation}\label{multiscale:y_cum_error}
  \|\hat{y}_{k,t} - \hat{y}_k\| \le \rho_f^{kT + t} \|\hat{e}_{0}\| + \frac{\alpha G_u^f(T-1)(G_u^J + L_h G_y^J)}{1 - \rho_f}.
  \end{equation}
  On the other hand, by recalling the definition of $e_k$ in Lemma \ref{lem:output_error}, we get
 \begin{align}\label{multiscale:h_error}
 \|h(u_{kT+t}) - \hat{y}_{k,t}\|  &\le \|h(u_{kT+t}) - y_{kT+t}\| + \|y_{kT+t} - \hat{y}_{k,t}\|\notag\\
 &= e_{kT+t} + \|y_{kT+t} - \hat{y}_{k,t}\|.
 \end{align}
 Substituting \eqref{multiscale:u_cum_error}, \eqref{multiscale:y_cum_error}, \eqref{multiscale:h_error} into \eqref{multiscale:gradient_h_error} and combining \eqref{multiscale:u_error}, we further derive
 \begin{align}
 &\quad\|\hat{u}_{k, t+1} - u_{kT + t + 1}\|\notag\\
 &\le (\sqrt{1 - 2\alpha\mu_J + \alpha^2L_{J,u}^2} + \alpha C_1)\|\hat{u}_{k,t} - u_{kT + t }\|\notag\\
 &\quad + \alpha C_2\|\hat{y}_{k,t} - y_{kT + t }\| + \alpha G_y^J C(e_{kT + t} + \rho_f^{kT} \|\hat{e}_{0}\|)\notag\\
 &\quad + \alpha^2G_{y}^JC(T-1)(G_u^J + L_h G_y^J)\left(1 + \frac{G_u^f}{1 - \rho_f}\right).\end{align}
 Consequently, we reach the linear matrix inequality
 \[\begin{bmatrix}\|\hat{u}_{k,t+1} - u_{kT + t + 1}\|\\
 \|\hat{y}_{k,t+1} - y_{kT + t + 1}\|\end{bmatrix}\le M\begin{bmatrix}\|\hat{u}_{k,t} - u_{kT + t }\|\\
 \|\hat{y}_{k,t} - y_{kT + t }\|\end{bmatrix} + \alpha G_y^JC\begin{bmatrix}\bar{e}_{kT}\\
 0\end{bmatrix} \]
 with
 \[\bar{e}_{kT}:=e_{kT + t} + \rho_f^{kT} \|\hat{e}_{0}\| + \alpha(T-1)(G_u^J + L_h G_y^J)\left(1 + \frac{G_u^f}{1 - \rho_f}\right).\]
 According to this recursive relation, we have
 \begin{equation}
 \|\hat{u}_{k+1} - u_{(k+1)T}\| \le \rho(M)^T \|\hat{u}_{k} - u_{kT}\| + \alpha G_y^JC\sum_{j=0}^{T-1}\rho(M)^j\bar{e}_{kT}.
 \end{equation}
 From the spectral radius condition, we get
 \begin{align}
 &\quad \limsup_{k \to \infty}\|u_{kT} - \hat{u}_k\|\notag\\
 &\le \frac{\alpha G_y^JC}{1-\rho(M)}\limsup_{k \to \infty}\bar{e}_{kT}\notag\\
 &\le \frac{\alpha^2 G_y^JC(G_u^J + L_h G_y^J)}{1-\rho(M)}\left(\frac{L_h}{1-\rho_f} + (T-1)\left(1 + \frac{G_u^f}{1 - \rho_f}\right)\right).
 \end{align}
 Finally, we reach \eqref{multiscale:error_bound} by Lemma \ref{lem:standard_fo}.
\end{appendix}


\bibliographystyle{plain}        
\bibliography{autosam}           



\end{document}